\documentclass[12pt]{article}
\date{}
\begin{document}
\font\aaa=msbm10 scaled \magstep 0
\title{Volume Growth and Curvature Decay  of
               Positively Curved K\"{a}hler Manifolds}
\author{Bing-Long Chen   and   Xi-Ping Zhu}
\maketitle
\begin{center}
     { Department of Mathematics,
 Zhongshan University
\vskip 1pt
 Guangzhou 510275, P. R. China
\vskip 1pt and
\vskip 1pt
      The Institute of Mathematical Sciences,
         The Chinese
\vskip 1pt
       University    of Hong Kong,  Hong Kong }
\end{center}
\begin{abstract}
 In this paper we obtain three results concerning the geometry of complete noncompact positively curved K\"{a}hler manifolds at infinity. The first
one  states that the order of volume growth of a complete noncompact
K\"{a}hler
manifold with positive bisectional curvature is  at least   half of
the real dimension (i.e., the complex dimension). The second one states
 that the curvature of a complete noncompact K\"{a}hler manifold
with positive bisectional curvature decays at least linearly in the
average sense. The third result is concerned with the relation between
the volume growth and the curvature decay. We prove that
the curvature decay of a complete noncompact K\"{a}hler manifold
with nonnegative curvature operator and with the maximal volume growth is precisely quadratic in certain average sense.
\end{abstract}
\vskip 0.6cm
\section{  Introduction}
\vskip 0.3cm
This paper is concerned with the geometric properties of positively
 curved complete noncompact manifolds near the infinities. Let us first consider the volume growth of the manifolds. When an $m$-dimensional complete noncompact Riemannian manifold has nonnegative Ricci curvature,
the classical Bishop volume comparison theorem implies that the volume
growth is at most as  the Euclidean volume growth. On the other hand,
Calabi and Yau (see [25]) showed that the volume growth of a complete noncompact $m$-dimensional Riemannian manifold with nonnegative Ricci curvature must be at least of linear, i.e.,
$$
Vol(B(x_0,r)) \geq cr, \mbox{\rm \quad  for all  \quad}  1 \leq r
<+ \infty,
$$
where $Vol(B(x_0,r))$ is the volume of the geodesic ball centered at
$x_0 \in M$ with radius $r$, and $c$ is some positive constant
depending on $x_0$ and the dimension $m$. The first result of this
paper is the following volume growth estimate for K\"{a}hler manifolds
with nonnegative holomorphic bisectional curvature.
\vskip 0.3cm
\noindent {\bf Theorem 1} {\it Let M be a complex n-dimensional complete
noncompact  K\"{a}hler manifold with nonnegative holomorphic bisectional
curvature. Suppose also its holomorphic bisectional curvature is
positive at least at one point. Then the volume growth of M satisfies
$$
Vol(B(x_0,r)) \geq cr^n, \quad \mbox{ for all} \quad 1 \leq r <
+\infty,
$$
where c is some positive constant depending on $x_0$ and the dimension
n.}
\vskip 0.3cm
Next we study the curvature behavior near the infinities.
We still consider complete noncompact K\"{a}hler manifolds with positive
holomorphic bisectional curvature. In view of the classical Bonnet-Myers
theorem, the Ricci curvature can not be uniformly bounded from below
by a positive constant. Our second result shows that the curvature
actually decays at least linearly in the average sense.
\vskip 0.3cm
\noindent {\bf Theorem 2}
{\it Let $M$ be a complete noncompact K\"{a}hler manifold with positive
holomorphic bisectional curvature. Then for any $x_0 \in M$, there exists a
positive constant $C$ such that
$$ \frac 1{Vol(B(x_0,r))} \int_{B(x_0,r)}R(x)dx \leq \frac C{1+r},
\quad \mbox{ for all } 0 \leq r < + \infty,$$ where $R(x)$ is the
scalar curvature of $M$.} \vskip 0.1cm The well-known conjecture
of Yau on uniformization theorems asks if a complete noncompact
K\"{a}hler manifold with  positive holomorphic bisectional
curvature is biholomorphic to a complex Euclidean space. The
combination of the above Theorem 2 and the main theorem of [8]
gives the following partial affirmative answer to the Yau's
conjecture. \vskip 0.3cm \noindent {\bf Corollary} {\it Let $M$ be
a complex two-dimensional complete noncompact K\"{a}hler manifold
with bounded and positive holomorphic bisectional curvature.
Suppose the volume growth of $M$ is maximal, i.e.,
$$Vol(B(x_0,r)) \geq cr^4, \mbox{  for all  } 0 \leq r < + \infty,$$
for some point $x_0 \in M$ and some positive constant $c$. Then $M$
is biholomorphic to the complex Euclidean space
$\aaa\mbox{C}^2$.}
\vskip 0.2cm
Further we consider the relation between the volume growth and the
 curvature decay.
In [18],
Mok-Siu-Yau studied complete K\"{a}hler manifold with nonnegative
holomorphic bisectional curvature and with the Euclidean volume growth.
 They showed that the curvatures can not decay faster than quadratically
unless the manifolds are flat. It was further shown in [9], [19]
this result is still valid without restriction on the volume
growth.
 On the other hand, Yau predicted in [26]
that if the volume growth is of Euclidean, then the curvature must
decay quadratically in  certain average sense. In [8], under an
additional assumption that the scalar curvature tends to zero at
infinity in the average sense, Tang and the authors confirmed this
for the complex surface case by using some special features in
dimension $2$ such as the Gauss-Bonnet-Chern formula and the
classification of holonomy algebras for  $4$-dimensional
Riemannian manifolds. The third result of the present paper is the
following affirmative answer to  this prediction of Yau for all
dimensions under a more restricted curvature assumption. \vskip
0.3cm \noindent {\bf Theorem 3}
   {\it Let M be a complex n-dimensional complete noncompact
K\"{a}hler manifold with bounded and nonnegative curvature operator.
Suppose that there exists a positive constant $c_1$ such that for a fixed
base point $x_0$, we have
$$
(*) \quad  Vol(B(x_0,r)) \geq c_1 r^{2n}, \quad \mbox{ for all}
\quad 0 \leq r < + \infty. \hskip 2.7cm
$$
 Then there exists a constant
$c_2 > 0$ such that  for all $x \in M$,  $r > 0$, we have
$$
\int_{B(x,r)} \frac {R(y)}{d(x,y)^{2n-2}} dy  \leq c_2 \log (2+r).
$$
In particular,
$$
\frac 1 {Vol(B(x_0,r))} \int_{B(x_0,r)}R(x)dx \leq c_3  \frac
{\log (2+r)}{r^2},  \quad \mbox{ for all} \quad 0 \leq r < + \infty,
$$
for some constant $c_3 > 0.$}
\vskip 0cm
There are plenty of complete noncompact convex surfaces in
$\aaa\mbox{R}^3$ with quadratic volume growth  (i.e., cone-like at infinity). All of them are  complex $1$-dimensional K\"{a}hler manifolds
satisfying the assumptions of Theorem 3. Clearly any product of finitely
many such K\"{a}hler manifolds still satisfies the assumptions of Theorem
3. In this way we can produce a lot of complete noncompact K\"{a}hler
manifolds satisfying all assumptions of Theorem 3.
\vskip 0cm
On a K\"{a}hler manifold, the curvature operator of the K\"{a}hler metric
reduces to the space ${\bigwedge}^{1,1}$ of $(1,1)$-forms, this is, the
curvature operator $R_m: {\bigwedge}^{1,1} \to {\bigwedge}^{1,1}$ and
$R_m \equiv 0$ on the subspace of ${\bigwedge}^2$ perpendicular to
${\bigwedge}^{1,1}$. Thus the best possibility for a K\"{a}hler metric
 is the positivity of its curvature operator restricted on the subspace
${\bigwedge}^{1,1}.$ In the last part of this paper we will
present a family of complete noncompact K\"{a}hler manifolds which
have strictly positive curvature operator when restricted on
$(1,1)$-forms and which satisfy the all assumptions of Theorem 3.
\vskip 0cm This paper contains seven sections. In Section 2 we
give the proof of Theorem 1. Section 3 is devoted to the proof of
Theorem 2. In Section 4
 we establish some basic estimates to the Ricci flow.
In Section 5 we obtain the crucial time  decay estimate for the Ricci
 flow and Theorem 3 will be proved in Section 6. In the last section we
will give a family of complete noncompact K\"{a}hler manifolds which have
strictly positive curvature operator restricted on $(1,1)$-forms, have
the Euclidean volume growth and have quadratic curvature decay.
\vskip 0cm
We are grateful to L.F.Tam for many helpful discussions and Professor
S.T.Yau for his encouragement.
\vskip 0.6cm
\section{ The Volume Growth Estimate}
\vskip 0.3cm
Let us first recall  a cut-off function constructed by Schoen and Yau
(see Theorem 1.4.2 in [20]).
\vskip 0.3cm
\newtheorem{Lemma}{Lemma}[section]
\setcounter{section}{2}
\begin{Lemma}  (Schoen-Yau [20])  Suppose M is an m-dimensional
complete Riemannian manifold with nonnegative Ricci curvature. Then there
exists a constant $C(m)>0,$ depending only on the dimension m, such that
for any $x_0 \in M$ and any number $0< r < + \infty,$ there exists a
smooth function $\varphi_r \in C^{\infty}(M)$ satisfying
$$
e^{-C(m)(1+ \frac {d(x,x_0)}{r})} \leq \varphi_r(x) \leq
e^{-(1+\frac {d(x,x_0)}{r})},
$$
$$
|\nabla \varphi_r (x)| \leq \frac {C(m)}r \varphi_r(x),
$$
$$
|\Delta \varphi_r(x)| \leq \frac {C(m)}{r^2} \varphi_r(x),
$$
for $x \in M$, where $d(x,x_0)$ is the distance between $x$ and $x_0$.
\end{Lemma}
\vskip 0.3cm
\noindent {\sc Proof.} \quad In [20], Schoen and Yau constructed a
$C^{\infty}$
function $f$ satisfying
$$
\frac 1{C}(1+d(x,x_0)) \leq f(x) \leq 1+d(x,x_0),
$$
$$
| \nabla f(x)| \leq C, \hskip 0.3cm
$$
$$
| \Delta f(x)| \leq C, \hskip 0.3cm
$$
for all $x \in M$ and some positive constant $C$ depending only on the
dimension $m$. Denote by $g_{ij}(x)$ the Riemannian metric of $M$. We
define a new metric on $M$ by
$$
\tilde{g}_{ij}(x)= \frac 1 {r^2}g_{ij}(x), \quad x \in M.
$$
Clearly the new metric $\tilde{g}_{ij}(x)$ is still a complete Riemannian
metric on $M$ with nonnegative Ricci curvature. Thus there exists a smooth
function $f_r(x) \in C^{\infty} (M)$ such that
$$
\frac 1 C(1+\frac {d(x,x_0)}{r}) \leq f_r(x) \leq 1+ \frac {d(x,x_0)}r,
$$
$$
| \tilde{\nabla} f_r(x)|_{\tilde{g}} \leq C,
$$
$$
|\tilde{\Delta} f_r(x) |_{\tilde{g}} \leq C,
$$
for all $x \in M$. Here $\tilde{\nabla}, \tilde{\Delta}$ and
$| \cdot |_{\tilde{g}}$ are the covariant derivative, Laplacian and norm with respect to the new metric $\tilde{g}_{ij}(x),$ and $d(x,x_0)$ is the geodesic distance
with respect to the original metric $g_{ij}(x).$ Therefore by setting
$\varphi_r (x) = \exp (-f_r(x))$ we get the desired cut-off function.
\hfill $\#$
\vskip 0.3cm
Let $M$ be a complex $n$-dimensional complete noncompact K\"{a}hler
manifold with nonnegative holomorphic bisectional curvature. Let
$x_0$ be an arbitrary point in $M$. Since $M$ is complete and noncompact,
there is a geodesic ray, say $\gamma (t),$ emanating from $x_0$. The
family of functions $\eta_t: M \to \aaa\mbox{R}$ defined by $\eta_t(x)
=d(x_0,\gamma(t))-d(x,\gamma(t))$ is Lipschitz continuous (with Lipschitz
constant 1) and also satisfies $| \eta_t(x)| \leq d(x,x_0)$
(by the triangle inequality). It is thus an equicontinuous family
uniformly
bounded on compact sets. By Ascoli-Arzela's theorem, a subsequence of
$\{ \eta_t \}$ converges to Lipschitz continuous (with Lipschitz
constant 1)
function $\eta : M \to \aaa\mbox{R},$ the convergence being uniform on
compact subsets of $M$. This function $\eta$ is called the Busemann
function
 of $\gamma$ (Cheeger-Gromoll [6]). It was shown in [24] (see Theorem
A(c) of [24]) that $\eta$ is plurisubharmonic and is strictly plurisubharmonic at the points where the  bisectional curvature
is positive.
We now use this plurisubharmonic function to prove Theorem 1.
\vskip 0.3cm
\noindent {\sc Proof of Theorem 1.} \quad Without loss of generality,
we may
assume that the holomorphic bisectional curvature is positive at $x_0$.
Then the Busemann function $\eta$ is strictly plurisubharmonic in a
neighborhood of $x_0$. Since the holomorphic bisectional curvature is
nonnegative, the Ricci curvature is also nonnegative. By Lemma 2.1 we
obtain the cut-off function $\varphi_r$ for any given $r>0$.
\vskip 0cm
Fix a small positive number $\delta$ and a large positive number $r$.
Let $\rho: \aaa\mbox{R}^{2n} \to \aaa\mbox{R}$ be a nonnegative smooth
function supported in the unit ball centered at the origin of
$\aaa\mbox{R}^{2n}$, with
$$\int_{\aaa\mbox{R}^{2n}} \rho(v) dv =1$$
We set
$$\rho_{\varepsilon}(v)=\frac 1 {{\varepsilon}^{2n}}\rho(\frac v{\varepsilon}),
\quad \mbox{\rm for } \quad v \in \aaa\mbox{R}^{2n}, \varepsilon >0,$$
\vskip 0cm
and
$$(\eta \ast \rho_{\varepsilon})(x) =\int_{\aaa\mbox{R}^{2n}}
\rho_{\varepsilon}(v) \eta (\exp_x (v))dv, \quad \mbox{\rm for } \quad
x \in M.$$
Clearly $\eta \ast \rho_{\varepsilon}$ is smooth and  converges to
$\eta$ uniformly on compact sets as $\varepsilon \to 0.$ It is well-known
(see for example [10]) that $\eta \ast \rho_{\varepsilon}$ is also
plurisubharmonic on any compact subset as $\varepsilon >0$ small enough.
Denote by $\omega$ the K\"{a}hler form of $M$. We compute
$$
\int\limits_{ \{ \varphi_r > \delta \} } (\varphi_r-\delta)^n(\sqrt{-1})^n
(\partial \bar{\partial}(\eta \ast \rho_{\varepsilon}))^n \hskip 6cm
$$
$$
\arraycolsep=1.5pt
\begin{array}[b]{rl}
= & \displaystyle -\int\limits_{ \{ \varphi_r > \delta \}} n(\varphi_r-\delta)^{n-1}
(\sqrt{-1})^n \partial \varphi_r \wedge \bar{\partial}(\eta
\ast \rho_{\varepsilon}) \wedge (\partial \bar{\partial}(\eta \ast
\rho_{\varepsilon}))^{n-1}
\\[4mm]
\leq & \displaystyle \int\limits_{ \{ \varphi_r > \delta \} }
\frac {2nC(2n)}r (\varphi_r- \delta)^{n-1} \varphi_r(\sqrt{-1})^{n-1}
(\partial \bar{\partial} (\eta \ast \rho_{\varepsilon}))^{n-1} \wedge
\omega
\\[4mm]
\leq & \displaystyle \int\limits_{ \{ \varphi_r > \delta \} }
\frac {(2nC(2n))^2}{r^2}(\varphi_r-\delta)^{n-2} \varphi_r^2
(\sqrt{-1})^{
n-2}(\partial \bar{\partial}(\eta \ast \rho_{\varepsilon}))^{n-2}
\wedge
\omega^2
\\[4mm]
\cdots
\\[4mm]
\leq & \displaystyle \int\limits_{ \{ \varphi_r > \delta \} }
\frac {(2nC(2n))^n}{r^n} \varphi_r^n \omega^n.
\end{array}
\eqno (2.1)
$$
Since $\eta$ is strictly plurisubharmonic in a neighborhood of $x_0$, by
letting $\varepsilon \to 0$ and then $\delta \to 0,$ we know that there
is a
positive number $c_0 >0,$ depending on $x_0$ and the dimension $n$, such that as $r \geq 1$,
$$
\arraycolsep=1.5pt
\begin{array}[b]{rl}
0 < c_0 \leq & \displaystyle \int_M \varphi_r^n(\sqrt{-1})^n(\partial
\bar{\partial}\eta )^n
\\[4mm]
\leq & \displaystyle \frac {(2nC(2n))^n}{r^n} \int_M \varphi_r^n
\omega^n
\\[4mm]
\leq & \displaystyle \frac {(2nC(2n))^n}{r^n}
\int_M \varphi_r \omega^n.
\end{array}
\eqno (2.2)
$$
By using the standard volume comparison and Lemma 2.1, we have
$$
\arraycolsep=1.5pt
\begin{array}[b]{rl}
\displaystyle \int_M \varphi_r \omega^n \leq & \displaystyle
\int_{B(x_0,r)}e^{-(1+\frac {d(x,x_0)}r)}
\omega^n
\\[4mm]
& \displaystyle +\sum\limits_{k=0}^{\infty} \int_{B(x_0,2^{k+1}r)
 \backslash B(x_0,2^kr)}
e^{-(1+\frac  {d(x,x_0)}r)}\omega^n
\\[4mm]
\leq & \displaystyle Vol(B(x_0,r))+ \sum\limits_{k=0}^{\infty}
e^{-2^k}(2^{k+1})^{2n} Vol(B(x_0,r))
\\[4mm]
\leq & \displaystyle CVol(B(x_0,r)),
\end{array}
\eqno (2.3)
$$
where $C$ is some positive constant depending only on the dimension $n$.
\vskip 0cm
Thus a combination of (2.2) and (2.3) gives the desired volume growth estimate. \hfill $\#$
\vskip 0.3cm
\noindent {\bf Remark 2.2}
Klembeck [15] (see also [23]) and Cao [3] (see Remark 1.3 in [3])
presented
some complete K\"{a}hler metrics on $\aaa\mbox{C}^n$  which have positive
holomorphic bisectional curvature everywhere such that the volume of the
geodesic ball $B(O,r)$ centered at the origin $O$ with respect to the
K\"{a}hler metric grows like $r^n$. Thus the volume growth estimate of
Theorem 1 is sharp.
\vskip 0.3cm
\noindent {\bf Remark 2.3}
The assumption that the holomorphic bisectional curvature is positive at least at one point is necessary. Indeed, let $M_1$  be a noncompact
convex surface in $\aaa\mbox{R}^3$ which is asymptotic to a cylinder at
infinity. Clearly $M_1$ is a complete noncompact Riemannian surface with
positive  curvature and has linear volume growth. Also let $\aaa\mbox{C}
P^{n-1}$ be the complex projective space with the Fubini-Study metric.
Then the product $M_1 \times \aaa\mbox{C}P^{n-1}$ is a complex $n$-dimensional complete noncompact K\"{a}hler manifold with nonnegative
holomorphic bisectional curvature. But its volume growth is of linear.
\vskip 0.3cm
\noindent {\bf Remark 2.4} In view of Theorem 1, it is naturally raised a question whether there is a similar volume growth estimate for Riemannian
manifold with positive sectional curvature. The answer is negative. More
precisely, for each dimension $n \geq 2$, there exist $n$-dimensional
 Riemannian manifolds which have positive sectional curvature everywhere
but
have linear volume growth. In fact on any bounded, convex and smooth
domain $\Omega$ in $\aaa\mbox{R}^n$, we can choose a strictly convex
 function $u(x)$ defined over $\Omega$ which tends to $+ \infty$ as $x$
approaches to the boundary $\partial \Omega$. The graph of the convex
function $u(x)$ is a hypersurface in $\aaa\mbox{R}^{n+1}$, denoted by
$M^n$. Clearly the hypersurface $M^n$ is strictly convex, i.e., the second
fundamental form is strictly positive definite. It then follows from the
Gauss equation that $M^n$ has strictly positive sectional curvature.
Since the domain $\Omega$ is bounded, the volume growth of
$M^n$ must be of linear.
\vskip 0.5cm
\section{The  Minimal Curvature Decay}
\vskip 0.3cm In this section, we will prove that the curvature of
complete positive curved K\"{a}hler manifold decays at least
linearly. \vskip 0.3cm \noindent {\sc Proof of Theorem 2.} Let
$\eta (x)$ be a Busemann function constructed from a ray
initiating from $x_0$. By a smoothing argument, we may assume
$\eta(x)$ is smooth, $| d \eta (x) | \leq 2$, and the Levi form
$\sqrt{-1} \partial \bar{\partial} \eta (x)$ of $\eta (x)$ is
strictly positive on $M$. \vskip 0.2cm Let $K_M$ be the canonical
line bundle of $M$, then there exists a continuous positive
function $\lambda (x)$ on $M$, such that for any $k>0$,
$$\sqrt{-1} \partial \bar{\partial} (k\eta (x)) +C(K_M)+Ric \geq
k\lambda (x)\omega(x) \eqno (3.1)$$ where $C(K_M)$ be the
curvature of $K_M$, and $\omega(x)$ is the K\"{a}hler form of $M$.
\vskip 0.2cm By the $L^2$-estimates of $\bar{\partial}$-operator
of H\"{o}rmander ([14]), there exist a positive constant $k_0 >
0$, and a nontrivial holomorphic section $S$ of $K_M$ such that
$$\int_M \| S \|^2e^{-k_0 \eta (x)} d V(x) < + \infty. \eqno (3.2)$$
and $\| S(x_0) \|^2=1.$
\vskip 0.2cm
Recall the Poincar$\acute{e}$-Lelong equation
$$\sqrt{-1} \partial \bar{\partial} \log \| S \|^2 =[S=0]
+Ric $$ in the sense of currents. And for any $0 < \delta < 1$, a
direct computation gives
$$\Delta \log (\| S \|^2 + \delta) \geq R(x)\frac {\|S\|^2}{\|S\|^2 + \delta}  \eqno (3.3)$$
and
$$\Delta \| S \|^2 \geq 0, \eqno
(3.4)$$ in the classical sense. By the mean-value-inequality (see
[17]), we have for any $x \in M,$ $r >0$,
$$\| S \|^2 (x) \leq \frac {c(n)}{Vol(B(x,r))}
\int_{B(x,r)} \| S \|^2.$$ Here and in the followings, $c(n)$ is
denoted by various positive constant depending only on the
dimension. \vskip 0.1cm Choosing $r=d(x,x_0)+1$ and noting that
$\eta(x) \leq 2(d(x,x_0)+1)$, we have
$$
\arraycolsep=1.5pt
\begin{array}{rl}
\| S \|^2 (x) \leq & \displaystyle \frac {c(n)}{Vol(B(x_0,1))}
\int_{B(x_0,1+d(x,x_0))} \| S \|^2
\\[4mm]
\leq & \displaystyle \frac {c(n)e^{2k_0d(x,x_0)}}{Vol(B(x_0,1))}
\int_M \| S \|^2 e^{-k_0 \eta}
\\[4mm]
 \leq & \displaystyle const.e^{2k_0 d (x,x_0)},
\end{array}
$$
hence
$$\log  \| S \|^2(x) \leq c_1 d(x,x_0)+c_2, \quad \mbox{ for some } c_1
, c_2 > 0. \eqno (3.5)$$ \vskip 0.1cm Suppose $M$ admits a
positive Green function $G=G(x,x_0)$. Then for any $\alpha >0,$
$\beta > 0$, $\varepsilon > 0$ and $0 < \delta < 1$, it follows
from (3.3) that
$$\displaystyle \int\limits_{ \{ \beta > G > \alpha \} } R(x)\frac
{\|S\|^2(x)}{\|S\|^2(x) + \delta}
(G(x,x_0)-\alpha)^{1+\varepsilon}
$$
$$
\arraycolsep=1.5pt
\begin{array}{rl}
 \leq & \displaystyle \int\limits_{ \{ \beta
> G > \alpha \} } \Delta \log (\| S \|^2 + \delta)(G- \alpha)^{1+
\varepsilon}
\\[4mm]
= & \displaystyle \int\limits_{ \{ \beta > G > \alpha \} } \log
(\| S \|^2 + \delta) \Delta (G-\alpha)^{1+ \varepsilon}
\\[4mm]
& \displaystyle \quad + \int\limits_{ \{ G= \beta \} } \frac
{\partial \log (\| S \|^2 + \delta)}{
\partial \vec{n}}(G-\alpha)^{1+\varepsilon}
\\[4mm]
& \displaystyle \quad -(1+ \varepsilon)\int\limits_{ \{ G= \beta
\} } \log (\| S \|^2 + \delta) (G-\alpha)^{\varepsilon} \frac
{\partial G}{\partial \vec{n}}
\end{array}
$$
where $\vec{n}$ is the outer normal vector of $\partial \{ \beta
> G > \alpha \}.$
\vskip 0.2cm Since $$\Delta
(G-\alpha)^{1+\varepsilon}=(1+\varepsilon)(G-\alpha)^{
\varepsilon-1} | \nabla G|^2 \geq 0, $$ on $\{\beta > G > \alpha
\}, $ we get
$$
\displaystyle \int_{ \{ \beta > G > \alpha \}} \log (\| S \|^2 +
\delta) \Delta (G-\alpha)^{1+\varepsilon}
$$
$$
\arraycolsep=1.5pt
\begin{array}{rl}
\leq & \displaystyle \sup\limits_{ \{ \beta
> G > \alpha \} } \log (\| S\|^2 + \delta) \int\limits_{ \{ \beta
> G > \alpha \}} \Delta(G-\alpha)^{1+\varepsilon}
\\[4mm]
= & \displaystyle \sup\limits_{ \{ \beta >G> \alpha \} } \log (\|
S \|^2 + \delta) \int\limits_{\{ G=\beta \} }(1+\varepsilon)
(G-\alpha)^{\varepsilon} \frac {\partial G}{\partial \vec{n}}.
\end{array}
$$
Letting $\varepsilon \to 0,$ it follows
$$
\arraycolsep=1.5pt
\begin{array}[b]{rl}
\displaystyle \int_{ \{ \beta > G > \alpha \}}R(x)\frac
{\|S\|^2(x)}{\|S\|^2(x) + \delta} (G(x,x_0)-\alpha) \leq &
\displaystyle \sup\limits_{ \{ G > \alpha \} } \log (\| S\|^2 +
\delta) \int\limits_{ \{ G=\beta \}} \frac {\partial G} {\partial
\vec{n}}
\\[6mm]
& \displaystyle +\int\limits_{ \{ G=\beta\}}(G-\alpha) \frac
{\partial \log (\|S\|^2 + \delta)}{\partial \vec{n}}
\\[6mm]
& \displaystyle -\int\limits_{\{ G=\beta\}} \frac {\partial
G}{\partial \vec{n}} \log (\| S \|^2 + \delta).
\end{array}
\eqno  (3.6)
$$
On $\{G= \beta \}$, $G$ and $ \frac {\partial G}{\partial \vec{n}}$ are
asymptotic to
$\frac {c(n)}{d^{2n-2}}$ and $(2n-2)c(n) \frac {1}{d^{2n-1}}$
respectively, as $\beta \to +\infty.$ So we have
$$| \int_{ \{ G=\beta \}}(G-\alpha)| \leq c(n) \frac {d^{2n-1}}{d^{2n-2}}
=c(n)d \to 0, \mbox{ as  } \beta \to +\infty,$$ and
$$\int_{\{ G=\beta \}}\frac {\partial G}{\partial \vec{n}} \to c(n),
 \mbox{ as } \beta \to + \infty.$$
Hence, letting $\beta \to +\infty $ in (3.6), it follows
$$\int_{ \{ G > \alpha \}} R(x)\frac {\|S\|^2(x)}{\|S\|^2(x) + \delta} (G(x,x_0)-\alpha) \leq
c(n) \sup\limits_{\{ G > \alpha \}} \log (\| S\|^2 + \delta)$$
since $\| S \|(x_0)=1.$ This implies
$$\int_{\{ G > 2 \alpha \} }R(x)\frac {\|S\|^2(x)}{\|S\|^2(x) + \delta}G(x,x_0) \leq c(n)
\sup\limits_{ \{ G > \alpha \}} \log (\|S \|^2 + \delta).$$ Then
letting $\delta \to 0$, we obtain, $$ \int_{\{ G > 2 \alpha \}
}R(x)G(x,x_0) \leq c(n) \sup\limits_{ \{ G
> \alpha \}} \log \|S \|^2. \eqno (3.7)$$ \vskip 0.2cm
Since the positive Green's function may not exist on $M$ in
general, we consider $\tilde{M}= M \times \aaa\mbox{C}^2,$
equipped with the product metric. By regarding $\|S\|^2$ as a
function on $\tilde{M}$, the inequality (3.3) still holds. Let
$z_0=(x_0,0),$ $z=(x,y) \in \tilde{M},$ the minimal positive Green
function of $\tilde{M}$ exists and satisfies
$$\frac {c(n)^{-1} \tilde{d}^2(z,z_0)}{Vol(\tilde{B}(z_0,\tilde{d}
(z,z_0))} \leq \tilde{G}(z,z_0) \leq \frac {c(n)
\tilde{d}^2(z,z_0)}{Vol(\tilde{B}(z_0,\tilde{d}(z,z_0))}, \eqno
(3.8)$$ where $\tilde{d}^2(z,z_0)=d^2(x,x_0)+|y|^2.$ \vskip 0.2cm
In fact, (3.8) comes from the following well-known estimate
([20]):
$$c(n)^{-1} \int_{\tilde{d}^2(z,z_0)}^{\infty} \frac {dt}{Vol((\tilde{B}
(z_0,\sqrt{t}))} \leq \tilde{G}(z,z_0) \leq c(n) \int_{\tilde{d}^2(z,z_0)}
^{\infty} \frac {dt}{Vol(\tilde{B}(z_0,\sqrt{t}))},$$
in particular, we have
$$\tilde{G}(z,z_0) \geq \frac {c(n)^{-1} a^2}{Vol(\tilde{B}(z_0,a))},
\quad \mbox{ for all } z \in \tilde{B}(z_0,a). \eqno (3.9)
$$
The same argument of  obtaining (3.7) gives us
$$\int\limits_{ \{ \tilde{G} > 2 \alpha \} }R(x) \tilde{G}(z,z_0)
\leq c(n) \sup\limits_{ \{ \tilde{G} > \alpha \}} \log \| S\|^2. \eqno
(3.10)$$
For $\alpha >0$, let $r (\alpha)$ be the maximum positive number such that
$$\tilde{B}(z_0,r(\alpha))  \subset \{
\tilde{G} > \alpha \}$$ which, together with (3.8), implies
$$\frac {c(n)^{-1}r^2(\alpha)}{Vol(\tilde{B}(z_0,r(\alpha)))}
< \alpha < \frac
{c(n)r^2(\alpha)}{Vol(\tilde{B}(z_0,r(\alpha)))}.$$
 Since
$$B(x_0,\frac 1 2 r(\alpha)) \times B(0,\frac 1 2 r(\alpha))
\subset \tilde{B}(z_0,r(\alpha)) \subset B(x_0,r(\alpha)) \times
B(0,r(\alpha)), \eqno (3.11)$$ we get by combining with (3.8) that
for any $z$ with $G(z,z_0)
> \alpha,$
$$\frac{c(n)^{-1}r^2(\alpha)}{Vol(\tilde{B}(z_0,r(\alpha)))}
\leq \frac
{\tilde{d}^2(z,z_0)}{Vol(\tilde{B}(z_0,\tilde{d}(z,z_0)))}$$ i.e.,
$$\frac {Vol(\tilde{B}(z_0,\tilde{d}(z,z_0)))}{Vol(\tilde{B}(z_0,r(\alpha)))}
\leq c(n) \frac {\tilde{d}^2(z,z_0)}{r^2(\alpha)},$$ and then
$$\frac {Vol(B(x_0,\frac 1 2 \tilde{d}(z,z_0)))}{Vol(B(x_0,r(\alpha)))}
\leq c(n) \frac {r^2(\alpha)}{\tilde{d}^2(z,z_0)}.$$ Thus if
$\tilde{d}(z,z_0) \geq 2 r (\alpha),$ we get that $r(\alpha) \geq
c^{-1}(n)d(z,z_0)$. So we get
$$\tilde{B}(z_0,r(\alpha)) \subset \{ \tilde{G}
> \alpha \} \subset {\tilde{B} (z_0,c(n)r(\alpha))}. \eqno (3.12)$$
By (3.9), (3.10), (3.11) and (3.12), we have
$$\frac {r^2(\alpha)}{Vol(B(x_0,r(\alpha)))}
\int_{B(x_0,r(\alpha))} R(x)dV(x) \leq c(n)
\sup\limits_{B(x_0,c(n)r(\alpha))} \log \| S\|^2.$$ Combining
(3.5), we get for any $a > 0$,
$$\frac 1{Vol(B(x_0,a))} \int_{B(x_0,a)}R(x)dV(x) \leq \frac C{a+1}.$$
Therefore the proof of the theorem is completed.

\vskip 0.5cm
\section{The Ricci Flow and Preliminary Estimates}
\vskip 0.3cm From now on we consider $(M, g_{_{\alpha
\bar{\beta}}})$ to be a complete noncompact K\"{a}hler manifold
satisfying the assumptions of Theorem 3. Our method to prove
Theorem 3 is not  directly working  on the  K\"{a}hler metric.
Instead we use the K\"{a}hler metric as initial data and evolve it
by Hamiltion's  Ricci flow. We study the dynamic property of the
evolving metric and prove that the curvature of the evolving
metric decays linearly in time. The curvature of the evolving
metric satisfies a nonlinear heat equation. Intuitively, the
Harnack inequality of heat equation bridges the time decay with
the space decay. By using the time decay estimate on the evolving
curvature we will be indeed able to prove that the curvature of
the initial metric decays quadratically in space in average sense.
In this section we give some preliminary estimates for the Ricci
flow. \vskip 0.1cm Let us  evolve the metric $g_{_{\alpha
\bar{\beta}}}$ according to the following Ricci flow equation
$$
\arraycolsep=1.5pt
\left\{
\begin{array}{rcl}
\displaystyle \frac {\partial g_{_{\alpha \bar{\beta}}}}{\partial t}(x,t) & = &
-R_{_{\alpha \bar{\beta}}}(x,t), \quad x \in M,t>0,
\\[4mm]
\displaystyle g_{_{\alpha \bar{\beta}}}(x,0) & = &g_{_{\alpha \bar{\beta}}}(x), \quad x \in M,
\end{array}
\right.
\eqno (4.1)
$$
where $R_{_{\alpha \bar{\beta}}}(x,t)$ denotes the Ricci curvature
tensor of the metric $g_{_{\alpha \bar{\beta}}}(x,t).$ \vskip 0cm
It was shown in [21] that the Ricci flow (4.1) has a maximal
solution $g_{\alpha \bar{\beta}}(x,t)$ on $M \times [0, T_{\max
})$ with either $T_{\max}=+ \infty$ or $0 < T_{\max} < + \infty$
and the curvature becomes unbounded as $t \to T_{\max} (< +
\infty).$ Since the initial K\"{a}hler metric has nonnegative
curvature operator, it is known from [11], [22] that the
nonnegativity of the curvature operator and K\"{a}hlerity are
preserved under the evolution of (4.1). \vskip 0cm When the
initial metric satisfies an additional assumption that the initial
curvature decays pointwisely to zero at infinity (i.e.,
$\lim\limits_{d(x,x_0) \to + \infty}R(x) =0$), one knows from
Hamilton [13] (Theorem 18.3 in [13]) that the Euclidean volume
growth condition is preserve under the Ricci flow. Thanks to
Theorem 2, we can remove the curvature decay assumption in the
following result. \setcounter{section}{4}
\begin{Lemma}
Suppose $(M,g_{_{\alpha \bar{\beta}}})$ is a complex
$n$-dimensional complete noncompact simply connected K\"{a}hler
manifold with bounded and nonnegative curvature operator and
satisfies the conditions (*)  in Theorem 3. Then the Euclidean
volume growth condition (*) is preserved under the evolution of
(4.1), i.e.,
$$
Vol_t(B_t(x,r)) \geq c_1 r^{2n}, \quad \mbox{for all} \quad
r>0, x \in M \quad and \quad t \in [0,T_{\max}),
$$
with the same constant $c_1$ in the condition (*). Here $B_t(x,r)$
is the geodesic ball of radius $r$ centered at $x$ with respect to
the metric $ g_{_{\alpha \bar{\beta}}}(\cdot,t),$ and the volume
$Vol_t$ is taken with respect to the metric $g_{_{\alpha
\bar{\beta}}}(\cdot,t).$
\end{Lemma}
\vskip 0.3cm
\noindent {\sc Proof.}  Let us first consider the evolution equation of the curvature operator of the evolving metric $g_{\alpha \bar{\beta}}(x,t)$.
By applying the strong maximum principle to the evolution equation, it was
shown by Hamilton [11] (see Theorem 8.3 of [11]) that there exists a small
positive constant $\delta_0( < T_{\max})$ such that on the time interval
$0 < t <\delta_0,$ the image of the curvature operator of
$g_{\alpha \bar{\beta}}(\cdot,t)$ is invariant under parallel translation
and constant in time. Since the curvature of $g_{\alpha \bar{\beta}}(\cdot,
\delta)$ is uniformly bounded for $\delta \in [0, \delta_0],$ it is easy to
see from the evolution equation (4.1) that each metric
$g_{\alpha \bar{\beta}}(\cdot,\delta)$ is quasi-isometric to the initial
metric $g_{\alpha \bar{\beta}}(\cdot)$ and satisfies
$$Vol_{\delta}(B_{\delta}(x_0,r)) \geq c_1 (\delta) r^{2n}, \quad
\mbox{ for all  }  0 \leq r < + \infty,  \delta \in [0,\delta_0],$$
where the function $c_1(\delta)$ is defined and positive
on $[0,\delta_0]$ such that $\lim\limits_{\delta \to 0}c_1(\delta)=c_1.$
Thus, without loss of generality, we may assume that the image of
the curvature operator of the initial metric $g_{\alpha \bar{\beta}}$
 is invariant under parallel translation. According to the decomposition
theorem of de Rham (see for example Theorem 8.1 in [16]), the
simply connected, complete K\"{a}hler manifold with the metric
$g_{\alpha \bar{\beta}}$ is  holomorphically isometric to the
direct product $\aaa\mbox{C}^k \times M_1 \times \cdots \times
M_m,$ where $M_1, \cdots, M_m$ are all simply connected, complete,
irreducible K\"{a}hler manifold. Since $M$ has nonnegative
curvature operator and satisfies the maximal volume growth
condition $(*)$, each $M_i$, $i=1, \cdots,m,$ also has nonnegative
curvature operator and is noncompact. By applying the well-known
theorem of Berger[1], each irreducible K\"{a}hler manifold $M_i$,
$i=1, \cdots, m,$ is either a Hermitian symmetric space or has its
holonomy group as $U(n_i),$ $SU(n_i)$ or $S_p(\frac {n_i}{2}),$
 where $n_i$ is the complex dimensional of $M_i$. Since $M_i$ has
nonnegative curvature operator, $M_i$ is thus either Ricci flat or
has $U(n_i)$ holonomy (see for example [5]). Hence each $M_i$,
$i=1, \cdots, m$ has $U(n_i)$ holonomy. By noting that $M_i$ is
simply connected, the restricted holonomy group agrees with the
full holonomy group. Also note that the image of the curvature
operator is invariant under parallel translation. We then have
from Ambrose-Singer holonomy theorem that the image of the
curvature operator at any point $x$ is $u(n_i) \cong
\bigwedge_{\aaa\mbox{R}}^{1,1}(T_xM_i),$ the space of real
$(1,1)$-forms on the tangent space $T_xM_i$. Therefore, together
with the nonnegativity of the curvature operator, it follows that
the curvature operator of $M_i$ is strictly positive when
restricted to $\bigwedge_{\aaa\mbox{R}}^{1,1}(T_xM_i).$ In
particular the holomorphic bisectional curvature of $M_i$ is
positive everywhere. \vskip 0.1cm Denote by $R_i(x)$ the scalar
curvature of the K\"{a}hler manifold $M_i$, $i=1, \cdots,m.$ We
then know from Theorem 2 that the scalar curvature $R_i(x)$ decays
at least linearly in the average sense. Since the scalar curvature
$R(x)$ of $M$ is given by $R_1(x)+ \cdots + R_m(x),$ it follows
that
$$\lim\limits_{r \to \infty} \frac 1{Vol(B(x_0,r))}
\int_{B(x_0,r)} R(x) dx=0, \eqno (4.2)$$
where $Vol(B(x_0,r))$ is the volume of the  geodesic ball
$B(x_0,r)$ with respect to the initial metric $g_{\alpha \bar{\beta}}.$
\vskip 0.2cm
  We have seen that the curvature operator of the evolving
metric $g_{_{\alpha \bar{\beta}}}(\cdot,t)$ is nonnegative, and then the
Ricci curvature $R_{_{\alpha \bar{\beta}}}(x,t)$ is also nonnegative on
$M \times [0, T_{\max})$. The equation in (4.1) thus implies that
the metric is shrinking in time. In particular,
$$g_{_{\alpha \bar{\beta}}}(x,t) \leq g_{_{\alpha \bar{\beta}}}(x,0),
\quad \mbox{\rm on}  \quad M \times [0, T_{\max}).$$
\vskip 0cm
Set
$$
F(x,t)= \log \frac {\det (g_{_{\alpha \bar{\beta}}}(x,t))}{\det
(g_{_{\alpha \bar{\beta}}}(x,0))}, \quad \mbox{\rm on } \quad M \times
[0,T_{\max}).
$$
We then have
$$
\arraycolsep=1.5pt
\begin{array}[b]{rl}
\displaystyle e^{F(x,t)}R(x,t)= & \displaystyle g^{_{\alpha \bar{\beta}}}(x,t)
R_{_{\alpha \bar{\beta}}}
(x,t) \cdot \frac {\det (g_{_{\alpha \bar{\beta}}}(x,t))}
{\det (g_{_{\alpha \bar{\beta}}}(x,0))}
\\[4mm]
\displaystyle \leq & \displaystyle g^{_{\alpha \bar{\beta}}}(x,0) R_{
_{\alpha \bar{\beta}}}(x,t)
\\[4mm]
\displaystyle = & \displaystyle g^{_{\alpha \bar{\beta}}}(x,0)(R_{_{\alpha \bar{\beta}}}(x,t)-
R_{_{\alpha \bar{\beta}}}(x,0))+R(x,0)
\\[4mm]
\displaystyle = & \displaystyle -\Delta_0 F(x,t)+R(x,0),
\end{array}
\eqno (4.3)
$$
where $\Delta_0$ denotes the Laplacian operator with respect to the initial
 metric $g_{_{\alpha \bar{\beta}}}(x,0)$ and $R(x,t)$ denotes the scalar
curvature of the metric $g_{_{\alpha \bar{\beta}}}(x,t).$ From the equation
in (4.1), we see
$$
\arraycolsep=1.5pt
\begin{array}[b]{rl}
\displaystyle \frac {\partial F(x,t)}{\partial t} = & \displaystyle
g^{_{\alpha \bar{\beta}}}(x,t) \frac {\partial}{\partial t}g_{_{\alpha
\bar{\beta}}}(x,t)
\\[4mm]
= & -R(x,t).
\end{array}
\eqno (4.4)
$$
A combination of (4.3) and (4.4) gives
$$
e^{F(x,t)} \frac {\partial F(x,t)}{\partial t}
\geq \Delta_0 F(x,t)-R(x,0), \quad {\rm on } \quad M \times [0,T_{\max}).
\eqno (4.5)
$$
\vskip 0cm
For any fixed point $x_0 \in M$ and  any number $0 < r < + \infty,$
we let $\varphi_r(x)$ be the cut-off function on $(M,g_{_{\alpha \bar{\beta}}})
$ obtained in Lemma 2.1. We compute
$$
\arraycolsep=1.5pt
\begin{array}{rl}
\displaystyle \frac {\partial}{\partial t} \int_M \varphi_r(x)
e^{F(x,t)} d V_0 \geq & \displaystyle \int_M \varphi_r(x) (\Delta_0
F(x,t)-R(x.0))d V_0
\\[4mm]
\geq & \displaystyle \frac {C(2n)}{r^2} \int_M \varphi_r(x) F(x,t)d V_0
-\int_M \varphi_r(x) R(x,0) dV_0,
\end{array}
$$
where $dV_0$ denotes the volume element of the initial metric $g_{_{\alpha
\bar{\beta}}}$ and we note that $F(x,t)$ is nonpositive. Note from (4.4)
that $F(\cdot,t)$ is nonincreasing in time and $F(\cdot,0) \equiv 0.$
We integrate the above inequality from $0$ to $t$ to get
$$
\int_M \varphi_r(x)(1-e^{F(x,t)})d V_0 \leq \frac {C(2n)t}{r^2}
\int_M \varphi_r(x)(-F(x,t))dV_0
$$
$$
\hskip 4cm +t \int_M \varphi_r(x) R(x,0)d V_0.
\eqno (4.6)
$$
Since the metric is shrinking under the Ricci flow, we have
$$
B_t(x_0,r) \supset B_0(x_0,r), \quad {\rm for } \quad t \geq 0, 0< r <
+ \infty,
$$
and
$$
\arraycolsep=1.5pt
\begin{array}[b]{rl}
\displaystyle Vol_t(B_t(x_0,r)) \geq & \displaystyle Vol_t(B_0(x_0,r))
\\[4mm]
= & \displaystyle \int_{B_0(x_0,r)}e^{F(x,t)} d V_0
\\[4mm]
= & \displaystyle Vol_0 (B_0(x_0,r))+\int_{B_0(x_0,r)}(e^{F(x,t)}-1)d V_0.
\end{array}
\eqno (4.7)
$$
By Lemma 2.1 and (4.6), the error term in (4.7) satisfies
$$
\arraycolsep=1.5pt
\begin{array}[b]{rl}
\displaystyle \int_{B(x_0,r)}(e^{F(x,t)}-1)dV_0 \geq & \displaystyle
e^{2C(2n)} \int_M \varphi_r(x) (e^{F(x,t)}-1)dV_0
\\[4mm]
\geq & \displaystyle \frac {C(2n)e^{2C(2n)}t}{r^2} \int_M \varphi_r(x)
F(x,t)dV_0
\\[6mm]
 & \displaystyle -e^{2C(2n)}t \int_M \varphi_r(x) R(x,0) dV_0.
\end{array}
\eqno (4.8)
$$
Consider any fixed $T_0 < T_{\max}.$ Since the curvature in uniformly bounded on $M \times [0,T_0],$ it is clear from (4.4) that $F(x,t)$ is
uniformly bounded on $M \times [0,T_0].$ We then set
$$
A= \sup \{ |F(x,t)| \quad | \quad x \in M, t \in [0,T_0] \}
$$
\vskip 0cm
and
$$
\varepsilon(r) = \sup \{ \frac 1{Vol_0(B_0(x_0,a))}
\int_{B_0(x_0,a)}R(x,0)dV_0  \quad | \quad a \geq r \} .
$$
From (4.2) we see   that $\varepsilon(r) \to 0$ as $r \to +\infty.$
By using the volume comparison theorem, we showed in (2.3) that
$$
\int_M \varphi_r(x)d V_0 \leq C Vol_0 (B(x_0,r)), \eqno (4.9)
$$
and similarly
$$
\arraycolsep=1.5pt
\begin{array}[b]{rl}
\displaystyle \int_M \varphi_r(x) R(x,0) dV_0 \leq & \displaystyle
\int_M R(x,0) e^{-(1+\frac {d_0(x_0,x)}r)}d V_0
\\[4mm]
\leq & \displaystyle \int_{B_0(x_0,r)}R(x,0)dV_0+\sum\limits_{k=0}^{\infty}
e^{-2^k}(2^{k+1})^{2n}
\\[6mm]
& \displaystyle
  \frac {Vol_0(B_0(x_0,r))}
{Vol_0(B_0(x_0,2^{k+1}r))}
\int_{B_0(x_0,2^{k+1}r)}R(x,0)d V_0
\\[6mm]
\leq & C \cdot \varepsilon(r) \cdot Vol_0 (B_0(x_0,r)),
\end{array}
\eqno (4.10)
$$
where $d_0(x_0,x)$ is the distance between $x_0$ and $x$ with respect
to the
initial metric $g_{_{\alpha \bar{\beta}}}(\cdot,0)$ and $C$ is some
positive
constant depending only on the dimension $n$.
\vskip 0.3cm
Substituting (4.8),(4.9) and (4.10) into (4.7) and dividing by $r^{2n},$
 we
obtain
$$
\arraycolsep=1.5pt
\begin{array}{rl}
\displaystyle \frac {Vol_t(B_t(x_0,r))}{r^{2n}} \geq & \displaystyle
\frac {Vol_0(B_0(x_0,r))}{r^{2n}}-\frac {C(2n)e^{2C(2n)}AT_0}{r^2}
(C \frac {Vol_0(B_0(x_0,r))}{r^{2n}})
\\[6mm]
& \displaystyle -e^{2C(2n)}T_0(C \cdot  \varepsilon(r) \cdot
\frac {Vol_0(B_0(x_0,r))}{r^{2n}})
\\[6mm]
= & \displaystyle  (1-\frac
{C(2n)e^{2C(2n)}AT_0}{r^2}-Ce^{2C(2n)}T_0 \varepsilon(r)) \frac
{Vol_0(B_0(x_0,r))}{r^{2n}}.
\end{array}
$$
Then by letting $r \to + \infty$, we deduce that
$$
\lim\limits_{r \to + \infty} \frac {Vol_t (B_t(x_0,r))}{r^{2n}} \geq c_1.
$$
Hence by using the standard volume comparison we get
$$
Vol_t(B_t(x,r)) \geq c_1r^{2n}, \quad \mbox{\rm for all } \quad
x \in M, 0 \leq r < +\infty  \quad {\rm and } \quad t \in [0,T_0].
$$
Finally,  \  since $T_0 \ (< T_{\max})$ is arbitrary,  \  this
completes the proof of the lemma.   \hfill $\#$ \vskip 0.3cm In
the next section we will use rescaling arguments to analyse the
behavior
 of the evolving metric  near the maximal time $T_{\max}$. In view of the
compactness theorem of Hamilton [12], we need to estimate the injectivity
radius of $(M,g_{_{\alpha \bar{\beta}}}(\cdot,t))$ in terms of the maximum
of the curvature.
\vskip 0.3cm
Let us recall the local injectivity radius estimate of Cheeger, Gromov and
Taylor [7] which says that for any complete Riemannian manifold $N$ of
dimension $m$ with $\lambda \leq \mbox{sectional curvature of } N \leq  \bigwedge$ and let $r$ be a positive constant satisfying
$r \leq   {\pi}/{4 \sqrt{\bigwedge}}$ if $\bigwedge >0$, then the injectivity
radius of $N$ at a point $x$ is bounded from below as follows
$$
inj(x,N) \geq r \frac {Vol(B(x,r))}{Vol(B(x,r))+V_{\lambda}^m(2r)},
$$
where $V_{\lambda}^m(2r)$ denotes the volume of a ball with radius $2r$
in the  $m$-dimensi-  \ onal space form $V_{\lambda}^m$ with constant sectional
curvature $\lambda.$
\vskip 0.2cm
Denote by
$$
R_{\max}(t) = \sup \{ R(x,t) \ | \ x \in M \}, \quad {\rm for} \quad
t \in [0,T_{\max}).
$$
By applying the local injectivity estimate to the evolving manifold
$(M,g_{_{\alpha \bar{\beta}}}
(\cdot,$ $  t)),$ we get
 $$
\arraycolsep=1.5pt
\begin{array}{rl}
inj(M,g_{_{\alpha \bar{\beta}}}(\cdot,t)) \geq & \displaystyle
\frac {\pi}{4 \sqrt{R_{\max}(t)}}
\frac {Vol_t(B_t(x,\frac {\pi}{4 \sqrt{R_{\max}(t)}}))}
{Vol_t(B_t(x,\frac {\pi}{4 \sqrt{R_{\max}(t)}}))+V_0^{2n}(1)
(\frac {\pi}{4 \sqrt{R_{\max}(t)}})^{2n}}
\\[6mm]
\geq & \displaystyle
\frac {\pi}{4 \sqrt{R_{\max}(t)}}(\frac {c_1}{c_1+V_0^{2n}(1)}),
\end{array}
$$
here we used Lemma 4.1. Thus we have proved
\vskip 0.2cm
\begin{Lemma}
Suppose $(M,g_{_{\alpha \bar{\beta}}})$ is assumed in Theorem 3 and let
$g_{_{\alpha \bar{\beta}}}(x, t),$ $ (x,t) \in M \times [0, T_{\max}),$ be the
maximal solution of (4.1). And suppose $M$ is simply connected. Then there exists a positive constant $\beta$
such that
$$
inj(M,g_{_{\alpha \bar{\beta}}}(\cdot,t)) \geq \frac {\beta}{ \sqrt{R_{\max}
(t)}},
$$
for $t \in [0, T_{\max}).$ \hfill $\#$
\end{Lemma}
\section{Curvature Decay in Time}
\vskip 0.3cm In this section we analyse the behavior of the
evolving metric near the maximal time $T_{\max}.$ \vskip 0.3cm
\newtheorem{Theorem}{Theorem}[section]
\setcounter{section}{5}
\begin{Theorem} Let $(M,g_{_{\alpha \bar{\beta}}}(x))$ be a complex
$n$-dimensional complete noncompact K\"{a}hler manifold satisfying the
assumptions of Theorem 3 and let $(M,g_{_{\alpha \bar{\beta}}}(x,t)),$
$t \in [0,T_{\max}),$ be the maximal solution of the Ricci flow (4.1)
with $g_{_{\alpha \bar{\beta}}}(x)$ as initial metric. Then
$T_{\max} = + \infty$ and the scalar curvature $R(x,t)$ of the solution
satisfies
$$
0 \leq R(x,t) \leq \frac C{1+t}, \quad { on} \quad M \times
[0,+ \infty)$$
for some positive constant C.
\end{Theorem}
\vskip 0.3cm \noindent {\sc Proof.} In [8] we established this
result for complex $2$-dimensional case under an additional
assumption that the initial
 curvature decays to zero at infinity, where we used some special features of complex $2$-dimension such as the
classification of holonomy algebras and the splitting results in the
(real)$4$-dimensional Riemannian manifolds obtained by the combination
of Berger [1] and Hamilton [11]. In the followings we don't assume the curvature decay and design a dimension
reduction procedure.
\vskip 0cm
According to Hamilton [13], the maximal solution of (4.1) is of either
one of the following types.
$$
\arraycolsep=1.5pt
\begin{array}{rl}
\mbox{Type I:} & \quad T_{\max} < + \infty \mbox{ and } \sup(T_{\max}-t)
R_{\max}(t) < +\infty;
\\[4mm]
\mbox{Type II:} & \quad \mbox{either } T_{\max} < +\infty
\mbox{ and } \sup (T_{\max}-t) R_{\max}(t)= +\infty,
\\[4mm]
& \quad \mbox{or } T_{\max}= +\infty \mbox{ and } \sup
t R_{\max} (t) = + \infty;
\\[4mm]
\mbox{Type III:} & \quad T_{\max}= + \infty \mbox{ and }
\sup t R_{\max} (t) < + \infty.
\end{array}
$$
\vskip 0cm
We need to show that the maximal solution must be of Type III. Let us
argue by contradiction. Suppose the maximal solution is of Type I or
Type II. We first note that we may assume $M$ is simply connected.
Indeed by considering the universal covering of $M$, the induced metric
of $g_{\alpha \bar{\beta}}(x,t)$ on the universal covering is clearly still a solution to the Ricci flow and satisfies all assumptions of
Theorem 3, and of course is still of Type I or Type II.In the previous section we have obtained an injectivity radius
estimate for the solution. By combining with a result of Hamilton (Theorem
16.4 and 16.5 in [13]), we know that there exists a sequence of dilations
of the solution converging to a limit $(\tilde{M}^{(1)}, \tilde{g}_{_{\alpha
\bar{\beta}}}^{_{(1)}}(x,t))$ which is a complete solution of the Ricci flow
with nonnegative curvature operator, exists for
$- \infty < t < \Omega$ for some $0 < \Omega \leq + \infty$ and
$$
\tilde{R}^{{(1)}}(x,t) \leq \Omega / (\Omega -t)
$$
holds everywhere with equality somewhere at $t=0$. Here we denote
by $\tilde{R}^{{(1)}}(x,t)$ the scalar curvature of the limiting
metric $\tilde{g}_{_{\alpha \bar{\beta}}}^{_{(1)}}(x,t),$ and
denote $\Omega / (\Omega-t)$ to be $1$ if $\Omega = + \infty.$
Clearly the limiting solution $(\tilde{M}^{(1)},
\tilde{g}_{_{\alpha \bar{\beta}}}^{_{(1)}}(x,t))$ is still
noncompact. By Lemma 4.1 and the rescaling invariance of the
condition (*), we see that the limiting solution
$(\tilde{M}^{(1)}, \tilde{g}_{_{\alpha
\bar{\beta}}}^{_{(1)}}(x,t))$ also satisfies
$$
Vol_t (\tilde{B}_t^{_{(1)}} (x,r)) \geq c_1 r^{2n} ,\quad \mbox{\rm
for all} \quad x \in \tilde{M}^{(1)} \quad \mbox{\rm and} \quad 0 \leq
r < +\infty, \eqno (5.1)
$$
where $Vol_t(\tilde{B}_t^{_{(1)}}(x,r))$ denotes the volume of the geodesic ball $\tilde{B}_t^{_{(1)}}(x,r)$ of radius $r$ centered at $x $ of
$\tilde{M}^{(1)}$ with respect to the metric $\tilde{g}_{_{\alpha
\bar{\beta}}}^{_{(1)}}(x,t).$
\vskip 0cm
Denote by $\tilde{d}_t^{_{(1)}}(x,x_0)$ the distance between two points
$x,x_0 \in \tilde{M}^{{(1)}}$ with respect to the metric $
\tilde{g}_{_{\alpha \bar{\beta}}}^{_{(1)}}(\cdot,t).$ We first claim that
at time $t=0$, we have
$$
\limsup\limits_{\tilde{d}_{_0}^{_{(1)}}(x,x_0) \to + \infty}
\tilde{R}^{{(1)}} (x,0) \tilde{d}_{0}^{_{(1)}2}   (x,x_0) = + \infty,
\eqno (5.2)
$$
for any fixed $x_0 \in \tilde{M}^{{(1)}}.$
\vskip 0.2cm
Suppose not, thus the curvature of the metric
$\tilde{g}_{_{\alpha \bar{\beta}}}^{_{(1)}} (\cdot,0)$ decays quadratically.
By applying a result of Shi (see Theorem 8.2 in [22]), we know that the
solution $\tilde{g}_{_{\alpha \bar{\beta}}}^{_{(1)}}(\cdot,t)$ exists for
all $t \in (- \infty, + \infty)$ and satisfies
$$
\lim\limits_{t \to \infty} \sup \{ \tilde{R}^{{(1)}}(x,t) \ | \ x
\in \tilde{M}^{(1)} \} =0.
$$
On the other hand, by the Li-Yau-Hamilton inequality of Cao [2], we have
$$
\frac {\partial \tilde{R}^{{(1)}}}{\partial t} \geq 0, \quad
{\rm on } \quad \tilde{M}^{(1)} \times (- \infty, +\infty).
$$
Thus we deduce that $\tilde{R}^{{(1)}} \equiv 0$ on $\tilde{M}^{(1)}
\times (-\infty ,+ \infty)$ which contradicts with the fact that the
scalar curvature $\tilde{R}^{{(1)}}$ achieves $1$ somewhere at $t =0$.
This proves the claim (5.2).
\vskip 0.2cm
With the estimate (5.2) we can then apply a lemma of Hamilton (Lemma 22.2
in [13]) to find a sequence of points $x_j, j=1,2, \cdots,$ in $\tilde{M}^{
(1)},$ a sequence of radii $r_j, j=1,2, \cdots,$ and a sequence of positive
numbers $\delta_j, j=1,2, \cdots,$ with $\delta_j \to 0$ such that
$$
\arraycolsep=1.5pt
\begin{array}{rl}
(a) & \quad \tilde{R}^{(1)}(x,0) \leq (1+\delta_j) \tilde{R}^{(1)}
(x_j,0), \mbox{ for } x \in \tilde{B}_0^{_{(1)}}(x_j,r_j);
\\[4mm]
(b) & \quad r_j^2 \tilde{R}^{(1)}(x_j,0) \to + \infty, \mbox{ as }
j \to + \infty;
\\[4mm]
(c) & \quad \mbox{if } s_j = \tilde{d}_{_0}^{_{(1)}}(x_j,x_0),
\mbox{ then } \lambda_j= s_j /r_j \to + \infty,
\mbox { as } j \to +\infty;
\\[4mm]
(d) &\quad \mbox{the balls } \tilde{B}_0^{_{(1)}}(x_j,r_j)
\mbox{ are disjoint.}
\end{array}
$$
Since the metric $\tilde{g}_{_{\alpha \bar{\beta}}}^{_{(1)}}(\cdot,0)$
has nonnegative curvature operator, the sectional curvature is
nonnegative. Denote the minimum of the sectional curvature of the metric
$\tilde{g}_{_{\alpha \bar{\beta}}}^{_{(1)}}(\cdot,0)$ at $x_j$ by
$\nu_j^{(1)}.$ We next claim that  the following holds
$$
\varepsilon_j = \frac {\nu_j^{(1)}}{ \tilde{R}^{(1)}(x_j,0) } \to 0, \mbox { as }
j \to + \infty. \eqno (5.3)
$$
In fact , suppose not, there exists a subsequence $j_k \to + \infty$
and some positive number $ \varepsilon >0$ such that
$$
\varepsilon_{j_k}= \frac {\nu_{j_k}^{(1)}}{\tilde{R}^{(1)}(x_{j_k},0)}
\geq \varepsilon, \mbox{ for all } k=1,2,\cdots. \eqno (5.4)
$$
We have seen that the scalar curvature $\tilde{R}^{(1)}(x,t)$ is pointwisely
 nondecreasing in time. Then by using the local derivative estimate of Shi
(see Theorem 13.1 in [13]) and (a), (b), we have
$$
\arraycolsep=1.5pt
\begin{array}[b]{rl}
\displaystyle \sup\limits_{x \in \tilde{B}^{(1)}(x_{j_k},r_{j_k})}
| \nabla \tilde{R}_m^{(1)}(x,0) |^2 \leq &
\displaystyle C(\tilde{R}^{(1)}(x_{j_k},0))^2(\frac 1{r_{j_k}^2}+
\tilde{R}^{(1)}(x_{j_k},0))
\\[4mm]
\leq & \displaystyle 2C(\tilde{R}^{(1)}(x_{j_k},0))^3,
\end{array}
\eqno (5.5)
$$
where $\tilde{R}_m^{(1)}$ is the Riemannian curvature tensor of
$\tilde{g}_{_{\alpha \bar{\beta}}}^{_{(1)}}$ and $C$ is a positive constant
depending only on the dimension.
\vskip 0cm
Denote by $\nu^{(1)}(x)$ the minimum of the sectional curvature of
$\tilde{g}_{_{\alpha \bar{\beta}}}^{_{(1)}}(\cdot,0)$ at $x$. From (5.4),
(5.5) and (b) we have
$$
\arraycolsep=1.5pt
\begin{array}{rl}
\nu^{(1)}(x) \geq & \displaystyle \nu_{j_k}^{(1)}-\sqrt{2C}(\tilde{R}^{
(1)}(x_{j_k},0))^{\frac 3 2 } \tilde{d}_0^{_{(1)}}(x,x_{j_k})
\\[4mm]
 \geq & \displaystyle \tilde{R}^{(1)}(x_{j_k},0)(\varepsilon-
\sqrt{2C} \cdot \sqrt{\tilde{R}^{(1)}(x_{j_k},0)} \cdot
\tilde{d}_0^{_{(1)}}(x,x_{j_k}))
\\[4mm]
\geq & \displaystyle \frac {\varepsilon}{2} \tilde{R}^{(1)}
(x_{j_k},0),
\end{array}
$$
as
$$
\tilde{d}_0^{_{(1)}}(x,x_{j_k}) \leq
\frac {\varepsilon}{2 \sqrt{2C} \cdot \sqrt{\tilde{R}^{(1)}}(x_{j_k},0)}
\hskip 5cm
$$
and $k$ large enough. Thus by combining with (a), there exists $k_0 >0$
such that for any $ k \geq k_0$ and $x \in
\tilde{B}_0^{(1)}(x_{j_k}, \frac {\varepsilon}{2 \sqrt{2C} \cdot \sqrt{
\tilde{R}^{(1)}(x_{j_k},0)}}),$ we have
$$
\frac {\varepsilon}{2} \tilde{R}^{(1)}(x_{j_k},0) \leq \mbox{ the
sectional curvature at } x \leq 2 \tilde{R}^{(1)} (x_{j_k},0).
$$
Therefore the balls $\tilde{B}_0^{(1)}(x_{j_k},\frac {\varepsilon}{
2 \sqrt{2C} \cdot \sqrt{\tilde{R}(x_{j_k},0)}}),$ $k_0 \leq k < + \infty,$
are a family of disjoint remote curvature $\beta-$bumps for some
$\beta >0$ in the sense of Hamilton [13]. But this contradicts with the
finite bumps theorem of Hamilton [13]. So we have proved the claim (5.3).
\vskip 0.3cm
\noindent In the followings we are going to rescale the limiting solution
$(\tilde{M}^{(1)},\tilde{g}_{_{\alpha \bar{\beta}}}^{_{(1)}}(x,t))$
along the points $x_j,j=1,2, \cdots.$ We begin with an injectivity
radius estimates. We claim that there exists a positive constant
$\alpha$ such that the injectivity radii of $(\tilde{M}^{(1)},
\tilde{g}_{_{\alpha \bar{\beta}}}^{_{(1)}}(\cdot,0))$ at $x_j$ $(j=
1,2,\cdots)$ are uniformly bounded from below by
$$
inj_{\tilde{M}^{(1)}}(x_j,\tilde{g}_{_{\alpha \bar{\beta}}}^{_{(1)}}
(\cdot,0)) \geq \frac {\alpha}{\sqrt{\tilde{R}^{(1)}(x_j,0)}},
\eqno (5.6)
$$
for $j= 1,2, \cdots.$
\vskip 0cm
We prove the claim by  contradiction. Suppose there exists a subsequence
$j_k$, $k=1,2, \cdots,$ such that
$$
\tilde{\varepsilon}_{j_k}= \sqrt{\tilde{R}^{(1)}(x_{j_k},0)}\cdot
inj_{\tilde{M}^{(1)}}(x_{j_k},\tilde{g}_{_{\alpha \bar{\beta}}}^{_{(1)}}
(\cdot,0)) \to 0, \eqno (5.7)
$$
as $k \to + \infty.$ For each $k$, let us choose $x_{j_k}$ as the new origin and dilate $(\tilde{M}^{(1)},\tilde{g}_{_{\alpha \bar{\beta}}}^{_{(1)}}
(\cdot,0))$  along $x_{j_k}$ such that the dilated manifold, says
$\tilde{M}_k^{(1)}$, has the injectivity radius $1$ at the origin and
has the curvature bounded by $2 \tilde{\varepsilon}_{j_k}^2$ on the ball
centered at the origin with the radius not less than
$\tilde{r}_{j_k}=r_{j_k} \sqrt{\tilde{R}^{(1)}(x_{j_k},0)}$ $(\to + \infty,$ by
(b)). We have seen that the  scalar curvature of $\tilde{g}_{_{\alpha
\bar{\beta}}}^{_{(1)}}(x,t)$ is pointwisely nondecreasing in time. The
curvature bounds on the balls at $t=0$ also give the bounds for previous
time in these balls. Then by the local derivative estimate of Shi (see Theorem
13.1 in [13]) and the convergence theorem of Hamilton [12] we  know that
a  subsequence of $\tilde{M}_k^{(1)},$ $k=1,2, \cdots,$ converges in $
C_{loc}^{\infty}$ topology to a complete noncompact flat manifold
which has the injectivity radius $1$ at the origin. But the estimate (5.1)
 and its rescaling invariance imply that the flat manifold is
$\aaa\mbox{R}^{2n}$ with the Euclidean metric. This contradictions  proves
 the claim (5.6).
\vskip 0cm
Now we let $x_j$ be the new origin, dilate the space by a factor
$\lambda_j$ so that the scalar curvature $\tilde{R}^{(1)}(x_j,0)$ becomes
$1$ at the origin at $t=0$, and dilate the time by
$\lambda_j^2$ so that it is still a solution to the Ricci flow.
The balls $\tilde{B}_0^{(1)}(x_j,r_j)$ are dilated to the balls centered at the origins of radii $\tilde{r}_j=r_j \sqrt{\tilde{R}^{(1)}(x_j,0)} \to
+ \infty$ as $j \to + \infty$ (by (b)). Since the scalar curvature
of the limit $\tilde{g}_{_{\alpha \bar{\beta}}}^{_{(1)}}(x,t)$ is
pointwise nondecreasing in time, the curvature bounds on
$\tilde{B}_0^{(1)}(x_j,r_j)$ also give bounds for previous time in these
balls. Thus  by applying the injectivity radius estimate (5.6) and the
compactness theorem of Hamilton [12] again we can get a limit for the
dilated solutions. It follows from (5.1), (5.3), (a) and (b) that the
limit is a complete noncompact solution, still denoted by
$( \tilde{M}^{(1)}, \tilde{g}_{_{\alpha \bar{\beta}}}^{_{(1)}}(x,t)),$
to the Ricci flow on $t \in ( - \infty,0]$ such that
$$
\arraycolsep=1.5pt
\begin{array}{rl}
(e)_1 & \quad \mbox{the curvature operator is still nonnegative;}
\\[1mm]
(f)_1 & \quad \tilde{R}^{(1)}(x,t) \leq 1, \mbox{ for }
x \in \tilde{M}^{(1)}, t \in (-\infty,0], \mbox{ and } \tilde{R}^{(1)}
(0,0)=1;
\\[1mm]
(g)_1 & \quad Vol_t( \tilde{B}_t^{(1)}(x,r)) \geq c_1r^{2n},
\mbox{ for all } x \in \tilde{M}^{(1)}, 0 \leq r <+ \infty;
\\[1mm]
(h)_1 & \quad \mbox{there exists a 2-plane at the origin so that at
} t=0, \mbox{ the correspond-}
\\[1mm]
& \quad \mbox{ing sectional curvature vanishes.}
\end{array}
$$
If we consider the universal covering of $\tilde{M}^{(1)}$, the
induced metric of $\tilde{g}_{_{\alpha
\bar{\beta}}}^{_{(1)}}(x,t)$ on the universal covering is clearly
still a solution to the Ricci flow and satisfies all of above
${\rm (e)_1, (f)_1, (g)_1,}$ and ${\rm (h)_1.}$ Thus, without loss
of generality, we may assume that $\tilde{M}^{(1)}$ is simply
connected. \vskip 0.1cm Since $\tilde{M}^{(1)}$ is a K\"{a}hler
manifold, at any point $x$, the curvature operator of the
underlying Riemannian manifold restricts to the real (1,1) forms,
this is, $\tilde{R}_m^{(1)}:
\bigwedge_{\aaa\mbox{R}}^{1,1}(T_x\tilde{M}^{(1)})
 \to \bigwedge_{\aaa\mbox{R}}^{1,1}(T_x\tilde{M}^{(1)})$
and $\tilde{R}_m^{(1)} \equiv 0$ on the subspace of $\bigwedge^2
(T_x\tilde{M}^{(1)})$ perpendicular to
$\bigwedge_{\aaa\mbox{R}}^{1,1} (T_x\tilde{M}^{(1)}).$ This is
equivalent to the fact that the Levi-Civita connection of the
underlying Riemannian metric on $\tilde{M}^{(1)}$ restricts to the
unitary frame bundle. If we define the Lie bracket on $\bigwedge^2
(T_x\tilde{M}^{(1)})$ as in [11] by
$$
[ \varphi, \psi]_{ij}=g^{kl} \varphi_{ik}\psi_{jl}-g^{kl}\psi_{ik}
\varphi_{jl}, \mbox{  where  } \varphi, \psi \in {\bigwedge}^2
(T_x\tilde{M}^{(1)}),
$$
then $\bigwedge_{\aaa\mbox{R}}^{1,1}(T_x\tilde{M}^{(1)})$ is a Lie
subalgebra of $\bigwedge^2(T_x\tilde{M}^{(1)})$ isomorphic to
$u(n)$. \vskip 0cm Let us consider the evolution  equation of the
curvature operator of
 $\tilde{g}_{_{\alpha \bar{\beta}}}^{_{(1)}}(x,$  $ t)$.  By applying the strong
maximum principle to the evolution equation as in [11] (see
Theorem 8.3 of [11]) we know that there exists a positive constant
$K$ such that on the time interval $- \infty < t < -K$, the image
of the curvature operator of $(\tilde{M}^{(1)},\tilde{g}_{_{\alpha
\bar{\beta}}}^{_{(1)}} (\cdot,t))$ at every point $x$ is  a fixed
Lie subalgebra of $\bigwedge_{\aaa\mbox{R}}^{1,1}
(T_x\tilde{M}^{(1)}) \cong u(n)$, invariant under parallel
translation and constant in time. \vskip 0cm We now want to show
that $\tilde{M}^{(1)}$ must be a reducible manifold. Suppose
$\tilde{M}^{(1)}$ is an irreducible manifold. Then according to
the well-known theorem of Berger [1], $\tilde{M}^{(1)}$ is either
a Hermitian symmetric space or has its holonomy group as $U(n),
SU(n)$ or $S_p(\frac n 2)$. Since $\tilde{M}^{(1)}$ has
nonnegative curvature operator, $\tilde{M}^{(1)}$ is thus either
Ricci flat or has $U(n)$ holonomy (see for example [5]). The case
of Ricci flat is ruled out by
 (e$)_1$ and  (f$)_1$. This says that $\tilde{M}^{(1)}$ has
$U(n)$ holonomy group. Since $\tilde{M}^{(1)}$ is simply
connected, the restricted holonomy group agrees with the full
holonomy group. We then have  from Ambrose-Singer holonomy theorem
that the image of the curvature operator $\tilde{R}_m^{(1)}$ at
every point $x \in \tilde{M}^{(1)}$ is $u(n) \cong \bigwedge_{
\aaa\mbox{R}}^{1,1}(T_x\tilde{M}^{(1)}).$ This, together with the
nonnegativity of $\tilde{R}_m^{(1)}$, implies that
$\tilde{R}_m^{(1)}$ is strictly positive when restricted to
$\bigwedge_{\aaa\mbox{R}}^{1,1} (T_x\tilde{M}^{(1)}).$ \vskip 0cm
Let $x$ be an arbitrary point in $\tilde{M}^{(1)}$ and $e_1,$
$e_2$ be two vectors at $x$. Denote by $J$ the complex structure
of $\tilde{M}^{(1)}$. Consider the following (complex) $(1,1)$
vector
$$(e_1+\sqrt{-1} Je_1) \wedge \overline{(e_2+\sqrt{-1} Je_2)} \hskip 3cm$$
$$
\arraycolsep=1.5pt
\begin{array}{rl}
\hskip 3cm
= & \displaystyle (e_1+\sqrt{-1}Je_1) \wedge (e_2- \sqrt{-1}Je_2)
\\[4mm]
= & \displaystyle (e_1 \wedge e_2 +J e_1 \wedge J e_2 ) +\sqrt{-1}
(Je_1 \wedge e_2- e_1 \wedge Je_2).
\end{array}
$$
Since $\tilde{R}_m^{(1)}$ is strictly positive when restricted to
$\bigwedge_{\aaa\mbox{R}}^{1,1}( \tilde{M}^{(1)}),$ we have
$$
\tilde{R}_m^{(1)}(e_1 \wedge e_2 +Je_1 \wedge Je_2,e_1 \wedge e_2+
J e_1 \wedge Je_2) >0.
$$
By the K\"{a}hlerity of $\tilde{M}^{(1)},$ the sectional curvature
$K(e_1,e_2)$ of the $2-$plane spanned by $e_1$, $e_2$ satisfies
$$
\arraycolsep=1.5pt
\begin{array}{rl}
K(e_1,e_2) = & \displaystyle \frac 1 4 \tilde{R}_m^{(1)} (e_1 \wedge e_2+
Je_1 \wedge Je_2, e_1 \wedge e_2 +Je_1 \wedge Je_2)
\\[4mm]
 > & 0.
\end{array}
$$
But the arbitrariness of the point $x$ and the vectors $e_1$,
$e_2$ gives a contradiction with (h)$_1$. So $\tilde{M}^{(1)}$ is
a simply connected reducible K\"{a}hler manifold. \vskip 0cm
According to the decomposition  theorem of de Rham (see for
example Theorem 8.1 in [16]) $\tilde{M}^{(1)}$ can be
isometrically splitted as the direct product $\tilde{M}_1^{(1)}
\times \tilde{M}_2^{(2)}$ where each fact is also K\"{a}hler.
Clearly each factor $\tilde{M}_i^{(1)},$ $i=1, 2,$ is still a
complete noncompact solution of the Ricci flow for $t \in ( -
\infty,0],$ has uniformly bounded and nonnegative curvature
operator and satisfies the associated Euclidean volume growth,
moreover at least one, denoted by $(\tilde{M}^{(2)},
\tilde{g}_{_{\alpha \bar{\beta}}}^{_{(2)}}(x,t)),$ $ t \in (
-\infty, 0],$ is not flat. Denote the complex dimension of
$\tilde{M}^{(2)}$ by $n_2$ which is strictly less than $n$. \vskip
0cm Now we repeat the above argument for $( \tilde{M}^{(2)},
\tilde{g}_{_{\alpha \bar{\beta}}}^{_{(2)}}(x,t)),$ $t \in (-
\infty,0].$ For clearity, we outline the main points as
followings. We first apply the Li-Yau-Hamilton inequality of Cao
[2] and a result of Shi (Theorem 8.2 in [22]) to show that for any
fixed $x_0 \in \tilde{M}^{(2)},$
$$
\limsup\limits_{\tilde{d}_0^{_{(2)}}(x,x_0) \to + \infty}
\tilde{R}^{(2)}(x,0) \tilde{d}_{0}^{_{(2)} 2}(x,x_0) = + \infty.
\eqno (5.2)^{\prime}
$$
By using a lemma of Hamilton (Lemma 22.2 in [13]) we can find a sequence
of points $x_j$, $j =1, 2, \cdots,$ in $\tilde{M}^{(2)},$ a sequence
of radii $r_j,$ $j=1,2, \cdots,$ and a sequence of positive number
$\delta_j,$ $j=1,2, \cdots,$ with $\delta_j \to 0$, such that the corresponding (a), (b), (c) and (d) hold. We then apply the finite
bumps theorem of Hamilton [13] to show that
$$ \frac {\nu_j^{(2)}}{\tilde{R}^{(2)}(x_j,0)} \to 0,
\mbox{ as } j \to + \infty, \eqno (5.3)^{\prime}
$$
where $\nu_j^{(2)}$ denotes the minimum of the sectional curvature of the
metric $\tilde{g}_{_{\alpha \bar{\beta}}}^{_{(2)}}(\cdot,0)$ at $x_j$.
Next by combining the property of Euclidean volume growth and a rescaling
argument we can show that the injectivity radii of $(\tilde{M}^{(2)},
\tilde{g}_{_{\alpha \bar{\beta}}}^{_{(2)}}(\cdot,$ $0))$ at $x_j$ are uniformly bounded from below by
$$
inj_{\tilde{M}^{(2)}}(x_j,\tilde{g}_{_{\alpha \bar{\beta}}}^{_{(2)}}
(\cdot,0)) \geq  \frac {\alpha}{\sqrt{\tilde{R}^{(2)}(x_j,0)}}, \eqno
(5.6)^{\prime}
$$
for $j=1, 2, \cdots.$ Here $\alpha$ is some positive constant. Thus we can
rescale the solution $(\tilde{M}^{(2)}, \tilde{g}_{_{\alpha \bar{\beta}}}^{
_{(2)}}(\cdot,t)),$ $t \in (-\infty,0],$ along the points $x_j$,
$j=1, 2, \cdots,$ to get a limit which is a complete noncompact solution, still
denoted by  $(\tilde{M}^{(2)},\tilde{g}_{_{\alpha,\bar{\beta}}}^{_{(2)}}
(x,t)),$ to the Ricci flow on $t \in (- \infty,0]$ such that
$$
\arraycolsep=1.5pt
\begin{array}{rl}
(e)_2 & \quad \mbox{the curvature operator is still nonnegative};
\\[2mm]
(f)_2 & \quad \tilde{R}^{(2)} (x,t) \leq 1, \mbox{ for } x \in
\tilde{M}^{(2)}, t \in (- \infty,0], \mbox{ and } \tilde{R}^{(2)}(0,0)=1;
\\[2mm]
(g)_2 & \quad \mbox{there exists some positive constant}\ c_2\
\mbox{such that}
\\[3mm]
& \qquad  Vol_t(\tilde{B}_t^{(2)}(x,r)) \geq c_2r^{2n_2}, \mbox{ for all }
x \in \tilde{M}^{(2)}, 0 \leq r < +\infty,
\\[3mm]
& \quad \mbox{where } \tilde{B}_t^{(2)}(x,r) \mbox{ is the geodesic ball of }\tilde{M}^{(2)} \mbox{ centered at } x \mbox{ and with}
\\[2mm]
& \quad \mbox{radius } r
\mbox{ with respect to the metric } \tilde{g}_{_{\alpha \bar{\beta}}}^{_{(2)}}
(\cdot,t);
\\[2mm]
(h)_2 & \quad \mbox{there exists a } \mbox{2-plane at the origin
so that at }t =0, \mbox{ the corresponding }
\\[2mm]
& \quad \mbox{sectional curvature vanishes.}
\end{array}
$$
Without loss of generality, we may assume that $\tilde{M}^{(2)}$ is simply
connected. Exactly as before by using the theorem of Berger [1] and the decomposition theorem of de Rham we further deduce that $\tilde{M}^{(2)}$
can be isometrically splitted as the direct product $\tilde{M}_1^{(2)}
\times \tilde{M}_2^{(2)}.$ Also each factor $\tilde{M}_i^{(2)},$
$i=1, 2, $ is still a complete noncompact solution of the Ricci flow for
$t \in (-\infty,0],$ has uniformly bounded and nonnegative curvature
operator and satisfies the associated Euclidean volume growth, moreover
at least one, denoted by $(\tilde{M}^{(3)}, \tilde{g}_{_{\alpha \bar{\beta}}}^{_{(3)}}(x,t)),$ $t \in (- \infty,0],$ is not flat. The complex
dimension of $\tilde{M}^{(3)}$ is strictly less than $n_2$.
\vskip 0cm
Hence by repeating these procedures we finally obtain a nonflat complex
$1-$dimensional complete noncompact solution, denoted by $(\tilde{M},
 \tilde{g}_{_{\alpha \bar{\beta}}}(x,t)),$ of the Ricci flow for
$t \in (- \infty,0]$, which has uniformly bounded and nonnegative curvature
and satisfies the following Euclidean volume growth
$$
Vol_t(\tilde{B}_t(x,r)) \geq \tilde{c} r^2, \quad \mbox{for all }
x \in \tilde{M}, 0 \leq r < +\infty, \eqno (5.8)
$$
where $\tilde{B}_t(x,r)$ denotes the geodesic ball of $\tilde{M}$
centered at $x$ of radius $r$ with respect to the metric
$ \tilde{g}_{_{\alpha \bar{\beta}}}
(\cdot,t),$ and $\tilde{c}$ is some positive constant. As the
curvature of $ \tilde{g}_{_{\alpha \bar{\beta}}}
 (x,t),$  is nonnegative, it follows from Cohn-Vossen inequality that
$$
\int_{\tilde{M}}\tilde{R}(x,t) d \sigma_t \leq 8 \pi , \eqno (5.9)
$$
where $\tilde{R}(x,t)$ is the scalar curvature of $(\tilde{M},
\tilde{g}_{_{\alpha \bar{\beta}}}(x,t))$  and $d \sigma_t$ is the
volume element of the metric $ \tilde{g}_{_{\alpha \bar{\beta}}}
(x,t).$ \vskip 0cm Now the metric $ \tilde{g}_{_{\alpha
\bar{\beta}}} (x,t))$  is a solution to the Ricci flow on the
Riemann surface $\tilde{M}$ over the ancient time interval $(-
\infty,0]$. Thus, (5.8) and (5.9) imply that for each $t < 0$, the
curvature of $ \tilde{g}_{_{\alpha \bar{\beta}}} (x,t)$  has
quadratic decay in the average sense of Shi [22], and then the a
priori estimate of Shi (see Theorem 8.2 in [22]) implies that the
solution $ \tilde{g}_{_{\alpha \bar{\beta}}} (x,t).$  exists for
all $t \in (- \infty , +\infty)$ and satisfies
$$
\lim\limits_{t \to + \infty} \sup \{ \tilde{R} (x,t) \ | \ x \in
\tilde{M} \} =0.
$$
Again by the Li-Yau-Hamilton inequality of Cao [2], we conclude that
$$
\tilde{R}(x,t) \equiv 0 , \mbox{ on } \tilde{M} \times (- \infty,
+ \infty).
$$
This contradicts with the fact that $ (\tilde{M},\tilde{g}_{_{\alpha
\bar{\beta}}} (x,t))$  is not flat.
\vskip 0cm
Therefore we have seeked the desired contradiction and have completed
the proof of Theorem 5.1. \hfill $\#$
\vskip 0.4cm
\section{Curvature Decay in Space}
\vskip 0.3cm
After obtain the time decay estimate, we can adapt the argument in [8]
to obtain the space decay estimates stated in Theorem 3. For sake of
completeness, we present a proof as follows.
\vskip 0.3cm
\noindent{\sc Proof of Theorem 3.}  Let $(M,g_{_{\alpha \bar{\beta}}})$
be  a complex $n-$dimensional complete noncompact K\"{a}hler manifold
satisfying all the assumptions of Theorem 3. Let
$ g_{_{\alpha \bar{\beta}}} (x,t)$ be the solution of the
Ricci flow (4.1) with $ g_{_{\alpha \bar{\beta}}} (x)$  as the initial
metric. From Theorem 5.1 we know that the solution $ g_{_{\alpha \bar{\beta}}} (x,t)$ exists
for all times $t \in [0, + \infty)$ and satisfies
$$
0 \leq R(x,t) \leq \frac {C}{1+t}, \quad
\mbox{on } M \times [0,+\infty), \eqno (6.1)
$$
for some positive constant $C$. Consider the following function
$$
F(x,t)=\log \frac {\det(g_{_{\alpha \bar{\beta}}} (x,t))}
{\det(g_{_{\alpha \bar{\beta}}} (x,0))} , \quad
x \in M ,t \geq 0,
$$
introduced in the proof of Lemma 4.1. Recall from the definition of
Ricci curvature tensor of a K\"{a}hler metric that
$$
\arraycolsep=1.5pt
\begin{array}{rl}
\displaystyle R_{\alpha \bar{\beta}}(x,t)-R_{\alpha \bar{\beta}}(x,0)=&
\displaystyle -\partial_{\alpha} \partial_{\bar{\beta}} \log \det (g_{\gamma
\bar{\delta}}(x,t))+\partial_{\alpha} \partial_{\bar{\beta}}
\log \det (g_{\gamma \bar{\delta}}(x,0))
\\[3mm]
= &  \displaystyle - \partial_{\alpha} \partial_{\bar{\beta}} F(x,t).
\end{array}
$$
Thus
$$
R(x,0) = \Delta_0 F(x,t) +g^{\alpha \bar{\beta}}(x,0)
R_{\alpha \bar{\beta}}(x,t), \quad x \in M,t \geq 0, \eqno (6.2)
$$
where $\Delta_0$ is the Laplacian operator of the metric
$g_{\alpha \bar{\beta}}(x,0).$
\vskip 0cm
Since $(M, g_{\alpha \bar{\beta}}(\cdot,0))$ has nonnegative Ricci
curvature and is of Euclidean volume growth, it is well known
(see for example [20]) that the Green function
$G_0(x,y)$ of the initial metric
$g_{\alpha \bar{\beta}}(\cdot,0)$ exists on $M$ and satisfies the estimates
$$
\frac {C_4^{-1}}{d_0^{2n-2}(x,y)} \leq G_0(x,y) \leq \frac
{C_4}{d_0^{2n-2}(x,y)}, \eqno (6.3)
$$
and
$$
| \nabla_y G_0(x,y)|_0 \leq \frac {C_4}{d_0^{2n-1}(x,y)}, \eqno (6.4)
$$
where $d_0(x,y),$ $| \cdot |_0$ are the geodesic distance of $x,$
$y$ and the norm with respect to the initial metric $g_{_{\alpha
\bar{\beta}}}(\cdot,0),$ and $C_4$ is some positive constant.
\vskip 0.2cm For any fixed $\bar{x}_0 \in M$ and any $\alpha > 0,$
we denote
$$
\Omega_{\alpha} = \{ x \in M \ | \ G_0 (\bar{x}_0,x) \geq \alpha \}.
$$
By (6.3), it is not hard to see
$$
B_0(\bar{x}_0, (\frac {C_4^{-1}}{\alpha})^{\frac 1 {2n-2}}) \subset
\Omega_{\alpha} \subset B_0(\bar{x}_0,(\frac {C_4}{\alpha})^{\frac
{1}{2n-2}}). \eqno (6.5)
$$
Recall from (4.4) that
$$
\frac {\partial F(x,t)}{\partial t}= - R(x,t), \quad x \in M, t \geq 0.
$$
By combining with the time decay estimate (6.1) we deduce that
$$
0 \geq F(x,t) \geq -C_5 \log (1+t), \quad x \in M, t \geq 0, \eqno (6.6)
$$
for some positive constant $C_5$.
\vskip 0.2cm
Multiplying (6.2) by $G_0(\bar{x}_0,x)-\alpha$ and integrating over
$\Omega_{\alpha},$ we have
$$ \int_{\Omega_{\alpha}} R(x,0)(G_0(\bar{x}_0, x)- \alpha)dx \hskip 5cm$$
$$\arraycolsep=1.5pt
\begin{array}[b]{rl}
= & \displaystyle \int_{\Omega_{\alpha}} (\Delta_0 F(x,t))(G_0(
\bar{x}_0,x)-\alpha)dx
\\[4mm]
& \quad \displaystyle +\int_{\Omega_{\alpha}} g^{\alpha \bar{\beta}}
(x,0) R_{\alpha \bar{\beta}}(x,t)(G(\bar{x}_0,x)-\alpha) dx
\\[4mm]
= & \displaystyle -\int_{\partial \Omega_{\alpha}}F(x,t) \frac
{\partial G_0 (\bar{x}_0,x)}{\partial \nu} d \sigma-F(\bar{x}_0,t)
\\[4mm]
& \quad \displaystyle +\int_{\Omega_{\alpha}} g^{\alpha \bar{\beta}}(x,0)
R_{\alpha
\bar{\beta}}(x,t)(G(\bar{x}_0,x)-\alpha)dx
\\[4mm]
\leq & \displaystyle C_5 (1+C_4^{1+\frac {2n-1}{2n-2}} \alpha^{\frac {2n-1}{2n-2}}
Vol_0(\partial \Omega_{\alpha})) \log (1+t)
\\[4mm]
& \quad \displaystyle + \int_{\Omega_{\alpha}}g^{\alpha \bar{\beta}}(x,0)
R_{\alpha \bar{\beta}}
(x,t) G_0(\bar{x}_0, x)dx,
\end{array}
\eqno(6.7)
$$
by (6.3), (6.4) and (6.6). Here we used $\nu$ to denote the outer unit
normal of $\partial \Omega_{\alpha}$. From the coarea formula, (6.3),
(6.4) and (6.5), we have
$$
\arraycolsep=1.5pt
\begin{array}{rl}
\displaystyle \frac 1 {\alpha} \int_{\alpha}^{2 \alpha}
r^{\frac {2n-1}{2n-2}} Vol_0( \partial \Omega_r)dr \leq
& \displaystyle 2^{\frac {2n-1}{2n-2}} \alpha^{\frac {1}{2n-2}}
\int_{\alpha}^{2 \alpha} \int_{\partial \Omega_r}
| \nabla G_0( \bar{x}_0,x)|_0 d \sigma |d \nu |
\\[4mm]
\leq & \displaystyle 2^{\frac {2n-1}{2n-2}} \cdot C_4^{1+{\frac {2n-1}
{2n-2}}} \cdot \alpha^{\frac {2n}{2n-2}}Vol_0(\Omega_{\alpha})
\\[4mm]
\leq
& \displaystyle 2^{\frac {2n-1}{2n-2}} \cdot C_4^{1+{\frac {2n-1}{2n-2}}}
\alpha^{\frac {2n}{2n-2}}Vol_0(B_0(\bar{x}_0,(\frac {C_4}{\alpha})^{\frac {1}{2n-2}})
\\[4mm]
\leq
& \displaystyle C_6
\end{array}
$$
for some positive constant $C_6$ by the standard volume comparison theorem.
Integrating (6.7) from $\alpha$ to $2 \alpha$ and using the above inequality,
we get
$$
\int_{\Omega_{\alpha}}R(x,0)(G_0(\bar{x}_0,x)-2 \alpha)dx
\leq C_5 (1+C_4^{1+{\frac {2n-1}{2n-2}}} \cdot C_6) \log (1+t) \hskip 3cm
$$
$$ \hskip 3cm +\int_{\Omega_{\alpha}}g^{\alpha \bar{\beta}}(x,0)
R_{\alpha
\bar{\beta}}(x,t)G_0(\bar{x}_0,x)dx. \eqno (6.8)
$$
It is easy to see that
$$
\int_{\Omega_{4 \alpha}}R(x,0)G_0(\bar{x}_0,x)dx \leq
2 \int_{\Omega_{2 \alpha}}R(x,0) (G_0(\bar{x}_0,x)-2\alpha)dx,
$$
and by the equation of the Ricci flow (4.1), we have
$$\int_0^t \int_{\Omega_{\alpha}}g^{\alpha \bar{\beta}}(x,0)
R_{\alpha \bar{\beta}}(x,t) G_0(\bar{x}_0,x)dxdt \hskip 3cm
$$
$$
\arraycolsep=1.5pt
\begin{array}{rl}
=& \displaystyle \int_{\Omega_{\alpha}}g^{\alpha \bar{\beta}}(x,0)
(g_{\alpha \bar{\beta}}(x,0)-g_{\alpha \bar{\beta}}(x,t))
G_0(\bar{x}_0,x)dx
\\[4mm]
\leq & \displaystyle n \int_{\Omega_{\alpha}}G_0(\bar{x}_0,x)dx.
\end{array}
$$
Thus by integrating (6.8) in time from $0$ to $t$ and combining the above
two inequalites, we get for any $t >0$,
$$ \int_{\Omega_{4 \alpha}}R(x,0) G_0(\bar{x}_0,x)dx
\leq 2 C_5(1+C_4^{1+{\frac {2n-1}{2n-2}}} \cdot C_6) \log (1+t)
+\frac {2n}{t} \int_{\Omega_{\alpha}}G_0(\bar{x}_0,x)dx.
$$
Finally, substituting (6.3) and (6.5) into the above inequality, we see
that
there exists some positive constant $C_7$ such that for any
$\bar{x}_0 \in M$, $t >0$ and $r >0$,
$$
\int_{B_0(\bar{x}_0,r)} \frac {R(x,0)}{d_0^{2n-2}(\bar{x}_0,x)}dx
\leq C_7 ( \log (1+t) +\frac {r^2}{t}).
$$
By choosing $t =r^2$ we get the desired first estimate. The second
estimate is a direct consequence of first estimate.
\hfill $\#$
\vskip 0.6cm
\section{Examples of Positively Curved K\"{a}hler }
\hskip 1.1cm {\Large\bf Manifolds} \vskip 0.3cm A homothetically
expanding gradient K\"{a}hler-Ricci soliton in a manifold $M$ is a
complete solution of the Ricci flow which moves along the equation
(4.1) by a one-parameter group of biholomorphisms in the direction
of a gradient holomorphic vector field and also expands by a
factor at the same time. More precisely, this means that the Ricci
tensor of $g_{\alpha \bar{\beta}}$ can be expressed  as
$$
\left\{
\arraycolsep=1.5pt
\begin{array}{rcl}
\displaystyle R_{\alpha \bar{\beta}} \quad & = & \displaystyle
\nabla_{\alpha} \nabla_{\bar{\beta}}f
-\rho g_{\alpha \bar{\beta}},
\\[4mm]
\displaystyle \nabla_{\alpha} \nabla_{\beta} f & = & \displaystyle
\nabla_{\bar{\alpha}} \nabla_{
\bar{\beta}} f=0,
\end{array}
\right.
\eqno (7.1)
$$
for some constant $\rho >0$ and some function $f$ on $M$. In [4], Cao
obtained a family of homothetically expanding gradient K\"{a}hler-Ricci
soliton on $\aaa\mbox{C}^n$ of positive sectional curvature by
constructing global potential functions on $\aaa\mbox{C}^n$. In this
section
we will further show that this family of K\"{a}hler-Ricci solitons
actually
have nonnegative curvature operator everywhere and have positive curvature
operator when restricted on the subspace of $(1,1)-$forms. Moreover the
soliton metrics have Euclidean volume growth and their curvatures are
quadratic decay.
\vskip 0.2cm
To begin we recall the construction of the K\"{a}hler-Ricci solitons on
$\aaa\mbox{C}^n$ in [4]. Let $z_1,z_2, \cdots,z_n$ denote the coordinate
functions on $\aaa\mbox{C}^n$ and denote $|z|= \sqrt{|z_1|^2+ \cdots+
|z_n|^2}$ the distance function from the origin $O$ with the Euclidean
 metric on $\aaa\mbox{C}^n.$ Let $t= \log |z|^2$ and also let
$u: \aaa\mbox{R} \to \aaa\mbox{R}$ be a smooth, convex and increasing
function. The K\"{a}hler potential function $u(t)$ on
$\aaa\mbox{C}^n$ induces a K\"{a}hler metric
$$
\arraycolsep=1.5pt
\begin{array}[b]{rl}
\displaystyle g_{\alpha \bar{\beta}} = & \displaystyle \partial_{\alpha}
\partial_{\bar{\beta}} u(t)
\\[4mm]
= & \displaystyle e^{-t} u^{\prime} \delta_{\alpha \beta}
+e^{-2t}z_{\beta} \bar{z}_{\alpha}(u^{\prime \prime}- u^{\prime}),
\end{array}
\eqno (7.2)
$$
on $\aaa\mbox{C}^n.$ It then follows easily that
$$
g^{\alpha \bar{\beta}} = e^t (u^{\prime} )^{-1} \delta_{\alpha \beta} +
\bar{z}_{\beta} z_{\alpha} ((u^{\prime \prime})^{-1} - (u^{\prime})^{-1})\eqno  (7.3)
$$
\vskip 0cm
and
$$
\det ( g_{\alpha \bar{\beta}})=e^{-nt} u^{\prime \prime}
(u^{\prime})^{n-1}.
\eqno (7.4)
$$
\vskip 0.2cm
Set
$$
\arraycolsep=1.5pt
\begin{array}{rl}
f(t) = & -\log \det (g_{\alpha \bar{\beta}})
\\[4mm]
= & nt -(n-1) \log u^{\prime}- \log u^{\prime \prime}.
\end{array}
$$
From the definition of Ricci curvature,
$$
R_{\alpha \bar{\beta}}+ g_{\alpha \bar{\beta}}=
\partial_{\alpha} \partial_{\bar{\beta}}(f+u).
\eqno (7.5)
$$
Consider the gradient vector field given by
$$
\arraycolsep=1.5pt
\begin{array}[b]{rl}
\displaystyle V^{\alpha}= & g^{\alpha \bar{\beta}}\partial_{\bar{\beta}}
(f+u)
\\[4mm]
= & \displaystyle (e^t(u^{\prime})^{-1} z_{\alpha}+e^tz_{\alpha}
((u^{\prime \prime})^{-1}-(u^{\prime})^{-1}))(f^{\prime}+u^{\prime})e^{-t}
\\[4mm]
= & \displaystyle \frac {f^{\prime}+ u^{\prime}}{u^{\prime \prime}}
z_{\alpha},
\end{array}
\eqno (7.6)
$$
which is holomorphic if and only if $\frac {f^{\prime}+ u^{\prime}}
{u^{\prime \prime}} $ is constant.
\vskip 0.1cm
Thus the  K\"{a}hler metric $g_{\alpha \bar{\beta}}= \partial_{\alpha}
\partial_{\bar{\beta}}u(t)$ is a solution of (7.1) with $\rho=1$ if
and only if
$$f^{\prime}+u^{\prime}= \lambda u^{\prime \prime}$$
for some constant $\lambda.$ Now by letting $u^{\prime} =\phi,$ the
equation (7.1) with $\rho=1$ for expanding K\"{a}hler-Ricci solitons is reduced to
$$\frac {\phi^{\prime \prime}}{\phi^{\prime}}+(\frac {n-1}{\phi}+ \lambda)
\phi^{\prime}=n+ \phi,
\eqno (7.7)
$$\vskip 0cm
i.e.,
$$(\log ( \phi^{\prime}
 \phi^{n-1}e^{\lambda \phi}))^{\prime} = n+ \phi,
$$
and then
$$
\phi^{\prime}=\frac 1{\lambda^{n+1} \phi^{n-1}}(\lambda^n \phi^n+
\sum\limits_{j=0}^{n-1}(-1)^{n-j} \frac {n!}{j!}(1- \lambda)(\lambda
\phi)^j +c e^{- \lambda \phi})
\eqno (7.8)
$$
where $c$ is a constant. It was shown in [4] that for each $\lambda >1$
and by choosing $c=(-1)^{n-1} n! (1-\lambda),$ the ODE (7.8) has a solution
with $\phi > 0$ and $ \phi^{\prime} >0$ such that the
K\"{a}hler metric $g_{\alpha \bar{\beta}}=\partial_{\alpha}
\partial_{
\bar{\beta}}u(t)$ is complete on $\aaa\mbox{C}^n$ and has positive
sectional curvature. In order to show these K\"{a}hler metrics
actually having positive curvature operator on the subspace of
$(1,1)-$forms, we compute the curvature  tensors as follows.
\vskip 0.1cm
By (7.2) and (7.3), a straightforward computation gives
$$
\arraycolsep=1.5pt
\begin{array}{rl}
\displaystyle \frac {\partial g_{_{\alpha \bar{\beta}}}}{\partial
z_{\gamma}}= & \displaystyle
e^{-2t}(u^{\prime \prime}-u^{\prime})(\bar{z}_{\alpha}
\delta_{\beta \gamma}+\bar{z}_{\gamma} \delta_{\alpha \beta})+
e^{-3t}(u^{(3)}-3u^{\prime \prime}+2u^{\prime}) \bar{z}_{\gamma}
z_{\beta} \bar{z}_{\alpha},
\\[4mm]
\displaystyle \frac {\partial^2 g_{_{\alpha \bar{\beta}}}}{\partial
z_{\gamma}
\partial \bar{z}_{\eta}}= & \displaystyle
e^{-2t}(u^{\prime \prime}-u^{\prime})(\delta_{\alpha \eta} \delta_{\beta
\gamma}+\delta_{\alpha \beta} \delta_{\gamma \eta})+
e^{-3t}(u^{(3)}-3u^{\prime \prime}+2u^{\prime}) \times
\\[4mm]
& \displaystyle [z_{\eta}\bar{z}_{\alpha} \delta_{\beta \gamma}+z_{\eta}
\bar{z}_{\gamma}
\delta_{\alpha \beta}+z_{\beta} \bar{z}_{\alpha}\delta_{\gamma \eta}
+\bar{z}_{\gamma} z_{\beta} \delta_{\alpha \eta}]
\\[4mm]
& \displaystyle +e^{-4t}(u^{(4)}-6u^{(3)}+11u^{\prime \prime}-6u^{\prime})
z_{\eta}
 \bar{z}_{\gamma}
z_{\beta} \bar{z}_{\alpha},
\end{array}
$$and
$$
\arraycolsep=1.5pt
\begin{array}{rl}
\displaystyle R_{\alpha \bar{\beta} \gamma \bar{\eta}}= & \displaystyle
-\frac {\partial^2 g_{_{\alpha \bar{\beta}}}}{\partial z_{\gamma}
\partial \bar{z}_{\eta}}+ g^{\xi \bar{\zeta}}
\frac {\partial g_{\xi \bar{\beta}}}{\partial \bar{z}_{\eta}}
\frac {\partial g_{\alpha \bar{\zeta}}}{\partial z_{\gamma}}
\\[4mm]
= & \displaystyle -e^{2t}(u^{\prime \prime}-u^{\prime})(
\delta_{\alpha \eta} \delta_{\beta \gamma}+ \delta_{\alpha \beta}
\delta_{\gamma \eta})
\\[4mm]
&\displaystyle  +e^{-3t}(-(u^{(3)}-3u^{\prime \prime}+2u^{\prime})
+(u^{\prime})^{-1}(u^{\prime \prime}-u^{\prime})^2) \times
\\[4mm]
& \displaystyle (\bar{z}_{\alpha}{z}_{\eta} \delta_{\beta \gamma}+z_{\eta}\bar{z}_{\gamma}
\delta_{\alpha \beta}+z_{\beta} \bar{z}_{\alpha}\delta_{\gamma \eta}
+\bar{z}_{\gamma} z_{\beta} \delta_{\alpha \eta})
\\[4mm]
& \displaystyle  +e^{-4t}[ (u^{\prime})^{-1}(u^{(3)}-3u^{\prime \prime}+
2u^{\prime})(u^{(3)}+u^{\prime \prime}-2u^{\prime})+
((u^{\prime \prime})^{-1}
\\[4mm]
& \displaystyle -(u^{\prime})^{-1})(u^{(3)}-u^{\prime \prime}
)^2
-(u^{(4)}-6u^{(3)}+11u^{\prime \prime}-6
u^{\prime})] \cdot \bar{z}_{\alpha} z_{\beta} \bar{z}_{\gamma}
z_{\eta}.
\end{array}
$$
\vskip 0.2cm
Note that the K\"{a}hler metric
$g_{\alpha \bar{\beta}}=\partial_{\alpha} \partial_{\bar{\beta}}u(t)$
is rotationally symmetric on $\aaa\mbox{C}^n.$ To facilitate computations,
we only need to compute the curvature at a point
$z=(r,0, \cdots, 0).$ Let $\xi^{\alpha \bar{\beta}}$ be any nonzero
$(1,1)-$vector at the point. From a straightforward computation we get
$$
\arraycolsep=1.5pt
\begin{array}{rl}
\displaystyle R_{\alpha \bar{\beta} \gamma \bar{\eta}} \xi^{\alpha
\bar{\beta}} \overline{\xi^{\gamma \bar{\eta}}}= & \displaystyle
e^{-2t}[a( | \sum\limits_{\alpha=1}^{n} \xi^{\alpha
\bar{\alpha}}|^2+ \sum\limits_{\alpha, \beta =1}^{n} |\xi^{\alpha
\bar{\beta}}|^2)
\\[4mm]
&\displaystyle +b(\sum\limits_{\alpha=1}^{n} |\xi^{1 \bar{\alpha}}|^2+
\overline{\xi^{1 \bar{1}}}(\sum\limits_{\alpha=1}^{n} \xi^{\alpha
\bar{\alpha}})+\xi^{1 \bar{1}} \overline{(\sum\limits_{\alpha=1}^{n}
\xi^{\alpha \bar{\alpha}})}
\\[4mm]
&\displaystyle + \sum\limits_{\alpha=1}^{n} |\xi^{\alpha
\bar{1}}|^2)+c| \xi^{1 \bar{1}}|^2]
\\[4mm]
= &\displaystyle e^{-2t}[(c+2a+4b) | \xi^{1 \bar{1}}|^2+(a+b)
(\xi^{1 \bar{1}} \overline{(\sum\limits_{\alpha=2}^{n} \xi^{\alpha
\bar{\alpha}})}
\\[4mm]
&\displaystyle +(\sum\limits_{\alpha=2}^{n} \xi^{\alpha \bar{\alpha}})
\overline{\xi^{1 \bar{1}}}+(a+b)(\sum\limits_{\alpha=2}^{n} |\xi^{1
\bar{\alpha}}|^2+\sum\limits_{\alpha=2}^{n} |\xi^{\alpha \bar{1}}|^2)
\\[4mm]
&\displaystyle +a(|\sum\limits_{\alpha=2}^{n} \xi^{\alpha \bar{\alpha}}|^2+
\sum\limits_{\alpha , \beta =2}^{n} |\xi^{\alpha \bar{\beta}}|^2),
\end{array}
$$
where the functions $a,$ $b $ and $c$ are defined by

$$
\arraycolsep=1.5pt
\begin{array}{rcl}
a & = & \phi -\phi^{\prime},
\\[4mm]
b & =& \phi^{-1}(\phi^{\prime}- \phi)^2-(\phi^{\prime \prime}
-3 \phi^{\prime}+2 \phi),
\\[4mm]
c & = & \phi^{-1}( \phi^{\prime \prime}-3 \phi^{\prime}+2 \phi)
( \phi^{\prime \prime}+\phi^{\prime}-2 \phi)+((\phi^{\prime})^{-1}
-\phi^{-1})(\phi^{\prime \prime}-\phi^{\prime})^2
\\[4mm]
& &-(\phi^{(3)}-6\phi^{\prime \prime}
+11 \phi^{\prime}-6 \phi).
\end{array}
$$
By Cauchy-Schwarz inequality, we see
$$
\arraycolsep=1.5pt
\begin{array}[b]{rl}
\displaystyle R_{\alpha \bar{\beta} \gamma \bar{\eta}} \xi^{\alpha
\bar{\beta}} \overline{\xi^{\gamma \bar{\eta}}} \geq
&\displaystyle e^{-2t}[(c+2a+4b) | \xi^{1 \bar{1}}|^2+(a+b)
(\xi^{1 \bar{1}} \overline{(\sum\limits_{\alpha=2}^{n} \xi^{\alpha
\bar{\alpha}})}
\\[4mm]
&\displaystyle +(\sum\limits_{\alpha=2}^{n} \xi^{\alpha \bar{\alpha}})
\overline{\xi^{1 \bar{1}}}+(a+b)(\sum\limits_{\alpha=2}^{n} |\xi^{1
\bar{\alpha}}|^2+\sum\limits_{\alpha=2}^{n} |\xi^{\alpha \bar{1}}|^2)
\\[4mm]
&\displaystyle +\frac n{n-1} a |\sum\limits_{\alpha=2}^{n} \xi^{\alpha
\bar{\alpha}}|^2],
\end{array}
\eqno (7.9)
$$
and we also see
$$a+b= \phi^{-1}((\phi^{\prime})^2- \phi \phi^{\prime \prime}),$$
$$c+2a+4b =(\phi^{\prime})^{-1}((\phi^{\prime \prime})^2-
\phi^{\prime} \phi^{(3)}).
$$
\vskip 0cm
Thus to show that the curvature operator is strictly positive on
the subspace $\bigwedge^{1,1}$ of $(1,1)$-forms we only need to verify
the following inequalities:
$$
\arraycolsep=1.5pt
\begin{array}{rl}
(A) & \quad a>0 ;
\\[2mm]
(B) & \quad a+b > 0;
\\[2mm]
(C) & \quad c+2a+4b > 0;
\\[2mm]
(D) & \quad \displaystyle (a+b)^2 < \frac {n}{n-1} a \cdot (c+2a+4b).
\end{array}
\hskip 5.5cm
$$
\vskip 0.3cm
\noindent {\sc Verification of (A).} \quad Recall from (7.8) with
$c=(-1)^{n+1} n! (1- \lambda)$ that
$$
\arraycolsep=1.5pt
\begin{array}[b]{rl}
\phi^{\prime} = & \displaystyle \frac 1 {\lambda^{n+1} \phi^{n-1}}
(\lambda^n \phi^n+ (-1)^{n+1} n! (1-\lambda) \sum\limits_{j=n}^{\infty}
\frac {(- \lambda \phi)^j}{j!})
\\[4mm]
= & \displaystyle \frac 1 {\lambda^{n+1} \phi^{n-1}}
(\lambda^{n+1} \phi^n+ (-1)^{n+1} n! (1-\lambda)
\sum\limits_{j=n+1}^{\infty}
\frac {(- \lambda \phi)^j}{j!}).
\end{array}
\eqno (7.10)
$$
We thus  have
$$
\arraycolsep=1.5pt
\begin{array}[b]{rl}
a = & \phi -\phi^{\prime}
\\[4mm]
= & \displaystyle \frac {(\lambda-1)n!}{\lambda^{n+1} \phi^{n-1}} \cdot
  (-1)^{n+1}
\sum\limits_{j=n+1}^{\infty}
\frac {(- \lambda \phi)^j}{j!} .
\end{array}
\eqno (7.11)
$$
Denote by
$$
f_k(x)=(-1)^k  \sum\limits_{j=k}^{\infty} \frac {(-x)^j}{j!},
\quad k=0,1, \cdots, n+1, \quad \mbox{ for } x \geq 0.
$$
Clearly
$$
f_{n+1}^{\prime} =f_n, \quad f_n^{\prime}=f_{n-1}, \quad \cdots,
\quad f_2^{\prime}=f_1, \quad f_1^{\prime}=f_0=e^{-x} >0.
$$
Then we deduce from (7.11) that
$$
a=\frac {(\lambda-1)n!}{\lambda^{n+1} \phi^{n-1}}f_{n+1}(
\lambda \phi) >0, \quad \mbox{as } \lambda > 1.
\eqno (7.12)
$$
\vskip 0.3cm
\noindent {\sc Verification of (B).} \quad By (7.7) and (7.11) we have
$$
\arraycolsep=1.5pt
\begin{array}{rl}
\displaystyle
\phi^{\prime \prime} \phi-( \phi^{\prime})^2= & \displaystyle
(n+\phi) \phi \phi^{\prime}-(n+\lambda \phi)(\phi^{\prime})^2
\\[4mm]
= & \displaystyle
\frac {(\lambda-1) \phi^{\prime}}{\lambda^{n+1} \phi^{n-1}}
[(n+\lambda \phi) n! f_{n+1}( \lambda \phi)-
(\lambda \phi)^{n+1}].
\end{array}
$$
Denote by
$$
g_k(x)=(k-1+x)(k-1)! f_k(x)-x^k, \quad k=1, \cdots, n+1, \quad
\mbox{ for } x \geq 0.
$$
It is easy to check that
$$
g_{n+1}^{\prime}=ng_n, \quad g_n^{\prime}=(n-1)g_{n-2}, \quad
\cdots, \quad g_2^{\prime}=g_1=-xe^{-x} \leq 0.
$$
Then we deduce that
$$
\arraycolsep=1.5pt
\begin{array}[b]{rl}
a+b= & \phi^{-1} (( \phi^{\prime})^2-\phi  \phi^{\prime \prime})
\\[4mm]
= & \displaystyle
\frac {(\lambda-1) \phi^{\prime}}{\lambda^{n+1} \phi^{n-1}}
(-g_{n+1}( \lambda \phi))
\\[4mm]
> & 0, \quad \mbox{ as } \lambda > 1.
\end{array}
\eqno (7.13)
$$
\vskip 0.3cm
\noindent {\sc Verification of (C).} \quad By (7.7) and (7.10) we have
$$
\arraycolsep=1.5pt
\begin{array}{rl}
\displaystyle \frac {\phi^{(3)} \phi^{\prime}-(\phi^{\prime \prime})^2}{(\phi^{\prime})^2}
= & \displaystyle (\frac {\phi^{\prime \prime}}{\phi^{\prime}})^{\prime}
\\[4mm]
= & \displaystyle (n+ \phi- (\frac {n-1}{\phi}+\lambda)
\phi^{\prime})^{\prime}
\\[4mm]
= & \displaystyle \frac {\phi^{\prime}}{\phi^2}[\phi^2+(n-1+(n-1+ \lambda
\phi)^2) \phi^{\prime}
\\[6mm]
& \displaystyle
-(n-1+\lambda \phi)(n+\phi) \phi]
\\[6mm]
= & \displaystyle \frac {\phi^{\prime}}{\phi^2}[\phi^2+(n-1+(n-1+ \lambda
\phi)^2)\frac 1{\lambda^{n+1} \phi^{n-1}}(\lambda^{n+1}\phi^n
\\[6mm]
& \displaystyle -(\lambda-1)n! f_{n+1}(\lambda \phi))-(n-1+\lambda \phi)
(n+\phi) \phi]
\\[6mm]
= & \displaystyle \frac {\phi^{\prime}(\lambda-1)}{(\lambda \phi)^{n+1}}
[(\lambda \phi)^{n+1}(n-2+ \lambda \phi)
\\[6mm]
& \displaystyle
-n!(n-1+(n-1+ \lambda \phi)^2)
 f_{n+1}(\lambda \phi) ].
\end{array}
$$
Note that
$$
f_n(x)=\frac {x^n}{n!}-f_{n+1}(x), \quad \mbox{for } x \geq 0.
$$
We thus have
$$
\arraycolsep=1.5pt
\begin{array}[b]{rl}
\displaystyle \frac {\phi^{(3)} \phi^{\prime}-(\phi^{\prime \prime})^2}{
(\phi^{\prime})^2}= & \displaystyle \frac {(\lambda-1) \phi^{\prime}}{(\lambda \phi)^{n+1}}[(n-1+(n-1+ \lambda \phi)^2)n! f_n( \lambda \phi)
\\[6mm]
& \displaystyle -n(n-1+\lambda \phi) (\lambda \phi)^n].
\end{array}
\eqno (7.14)
$$
Denote by
$$h_k(x)=(k-1+(k-1+x)^2)k!f_k(x)-k(k-1+x)x^n,$$
$$\qquad k =1,2 ,\cdots,n, \quad \mbox{ for } x \geq 0.$$
It is easy to check that
$$h_n^{\prime}=nh_{n-1}, \quad h_{n-1}^{\prime}=
(n-1)h_{n-2} ,\quad \cdots, \quad h_2^{\prime}=2h_1=2x^2(-e^{-x}) \leq 0.
$$
Thus we deduce from (7.14) that
$$
\arraycolsep=1.5pt
\begin{array}[b]{rl}
c+2a+4b = & (\phi^{\prime})^{-1}((\phi^{\prime \prime})^2-
\phi^{\prime} \phi^{(3)})
\\[4mm]
=& \displaystyle \frac {(\lambda-1) (\phi^{\prime})^2}{(\lambda \phi)^{n+1}}(-h_n( \lambda \phi))
\\[4mm]
> & 0, \qquad \mbox{ as } \lambda > 1.
\end{array}
\eqno (7.15)
$$
\vskip 0.3cm
\noindent {\sc Verification of (D).} \quad From (7.12), (7.13) and (7.15),
we have
$$\frac n{n-1} a \cdot (c+2a+4b)-(a+b)^2 \hskip 7.5cm$$
$$
\arraycolsep=1.5pt
\begin{array}{rl}
= & \displaystyle \frac {(\phi^{\prime})^2(\lambda-1)^2}{(\lambda^{n+1} \phi^n)^2}
[\frac n{n-1} n! f_{n+1}(\lambda \phi)(n(n-1+\lambda \phi)(
\lambda \phi)^n
\\[4mm]
& \quad \displaystyle -(n-1+(n-1+\lambda \phi)^2)n!f_n(\lambda \phi))
\\[4mm]
& \quad \displaystyle -((\lambda \phi)^{n+1}-(n+\lambda \phi) n!f_{n+1}(\lambda \phi))^2]
\\[4mm]
= & \displaystyle \frac {(\phi^{\prime})^2(\lambda-1)^2}{(\lambda^{n+1} \phi^n)^2}
[\frac n{n-1} n! (\frac {(\lambda \phi)^n}{n!}-f_{n}(\lambda \phi))(n(n-1+\lambda \phi)(
\lambda \phi)^n
\\[4mm]
& \quad \displaystyle -(n-1+(n-1+\lambda \phi)^2)n!f_n(\lambda \phi))
\\[4mm]
& \quad \displaystyle -((\lambda \phi)^{n+1}-(n+\lambda \phi) n!
(\frac {(\lambda \phi)^n}{n!}-f_{n}(\lambda \phi))^2]
\\[4mm]
= & \displaystyle \frac {(\phi^{\prime})^2(\lambda-1)^2}{(\lambda^{n+1} \phi^n)^2}
[\frac n{n-1}  (n(n-1+\lambda \phi)(
\lambda \phi)^{2n}
\\[4mm]
& \quad \displaystyle -((n-1+(n-1+\lambda \phi)^2)+n(n-1+\lambda \phi))(\lambda \phi)^n
n!f_n(\lambda \phi)
\\[4mm]
& \quad \displaystyle +(n!)^2(f_{n}(\lambda \phi))^2(n-1+(n-1+\lambda
 \phi)^2))-(n+ \lambda \phi)^2(n!)^2(f_{n}(\lambda \phi))^2
\\[4mm]
& \quad \displaystyle -2(n+ \lambda \phi)n!f_{n}(\lambda \phi) \cdot
n( \lambda \phi)^n + n^2(\lambda \phi)^{2n}]
\\[4mm]
= & \displaystyle \frac {n^2(\phi^{\prime})^2(\lambda-1)^2}{(n-1)
(\lambda^{n+1} \phi^n)^2}[ ((n-1)! f_{n}(\lambda \phi))^2(n(n-1
\\[4mm]
& \quad \displaystyle +(n-1+\lambda \phi)^2)-(n-1)(n+
\lambda \phi)^2)
\end{array}
$$
$$
\arraycolsep=1.5pt
\begin{array}[b]{rl}
& \quad \displaystyle  -(n-1)!f_n(\lambda \phi)
((n-1+(n-1+\lambda \phi)^2
\\[4mm]
& \quad \displaystyle +n(n-1+\lambda \phi))(\lambda \phi)^n
-2(n-1)(n+\lambda \phi)(\lambda \phi)^n)
+(\lambda \phi)^{2n+1}]
\\[4mm]
= & \displaystyle \frac {n^2(\phi^{\prime})^2(\lambda-1)^2(f_{n}(\lambda \phi))^2
}{(n-1)
(\lambda^{n+1} \phi^n)^2}[ ((n-1)! f_{n}(\lambda \phi))^2
\\[6mm]
& \quad  \displaystyle -(n-1)! f_{n}(\lambda \phi)((\lambda \phi)^n
+n(\lambda \phi)^{n-1})+(\lambda \phi)^{2n-1}].
\end{array}
\eqno (7.16)
$$
Thus to verify (D) it is sufficient to prove the following inequality
$$
((n-1)!f_n(x))^2-(n-1)!f_n(x)(x^n+nx^{n-1})+x^{2n-1} > 0, \eqno
(7.17)
$$
as  $x > 0.$
\vskip 0.2cm
Set
$$
L(x) =((n-1)!e^xf_n(x))^2-e^x(x^n+nx^{n-1})((n-1)!e^xf_n(x))
+e^{2x}x^{2n-1},
$$
for $x \geq 0.$ By noting that
$$
\arraycolsep=1.5pt
\begin{array}{rl}
\displaystyle (e^x f_n(x))^{\prime}= & \displaystyle
e^x(f_n(x)+f_{n-1}(x))
\\[4mm]
= & \displaystyle e^x \cdot \frac {x^{n-1}}{(n-1)!}
\end{array}
$$
We compute
$$
\arraycolsep=1.5pt
\begin{array}[b]{rl}
\displaystyle (L(x))^{\prime}= & \displaystyle 2e^{2x}((n-1)!)^2f_n(x) \cdot \frac {x^{n-1}}{(n-1)!}
-e^{2x}(x^n+nx^{n-1})(n-1)!\frac {x^{n-1}}{(n-1)!}
\\[4mm]
&\displaystyle -e^{2x}(x^n+2nx^{n-1}+n(n-1)x^{n-2})(n-1)!f_n(x)
\\[4mm]
&
\displaystyle +
2e^{2x}x^{2n-1}+(2n-1)e^{2x}x^{2n-2}
\\[4mm]
= & \displaystyle e^{2x}x^{n-2}[(n-1)!f_n(x)((2-2n)x-x^2-n(n-1))
\\[4mm]
& \displaystyle
+
x^{n+1}+(n-1)x^n].
\end{array}
\eqno (7.18)
$$
Denote by
$$
l(x)=(n-1)!f_n(x)((2-2n)x-x^2-n(n-1))+x^{n+1}+(n-1)x^{n},
 \mbox{ for } x \geq 0.
$$
By a direct computation,
$$
\arraycolsep=1.5pt
\begin{array}{rl}
\displaystyle (l(x))^{(n)}= & \displaystyle -(n-1)! \sum\limits_{k=0}^{n}
{n \choose k}
(f_n(x))^{(n-k)}(2(n-1)x+x^2+n(n-1))^{(k)}
\\[4mm]
& \displaystyle +(n+1)!x+(n-1) \cdot n!
\\[4mm]
= &\displaystyle -(n-1)![e^{-x}(2(n-1)x+x^2+n(n-1))+n(1-e^{-x})(2(n-1)
\\[4mm]
& \displaystyle +
2x)
+n(n-1)(e^{-x}-1+x)]+(n+1)! x +(n-1) \cdot n!
\\[4mm]
=& \displaystyle (n-1)! e^{-x}(2x-x^2).
\end{array}
$$
Thus
$$
(l(x))^{(n-1)}=(n-1)!x^2 e^{-x} >0, \quad \mbox{ for } x >0.
$$
Noting that
$$
l(0)=l^{\prime}(0)= \cdots=l^{(n-1)}(0)=0,
$$
we deduce that
$$
l(x) > 0, \quad \mbox{ for } x>0.
$$
Hence by combining with (7.18) we know that
$$
\arraycolsep=1.5pt
\begin{array}{rl}
L(x)= & \displaystyle e^{2x}[((n-1)!f_n(x))^2-(n-1)!f_n(x)
(x^n+nx^{n-1})+x^{2n-1}]
\\[4mm]
 > & 0, \qquad \mbox{ for } x>0.
\end{array}
$$
So we have verified the inequality (D).
\vskip 0.3cm
Summarizing above,   \ we have shown that the complete K\"{a}hler metrics
$g_{\alpha \bar{\beta}}=\partial_{\alpha}\partial_{\bar{\beta}}u(t)$
obtained by solving (7.8) for each $\lambda > 1$ with
$c=(-1)^{n+1}n!(1- \lambda)$ have strictly positive curvature operator
on the subspace of $(1,1)-$forms.
\vskip 0.2cm
Finally we study the volume growth and the curvature decay of the K\"{a}hler
metrics. From the equation (7.8) one can see that $\phi$ tends to
$+ \infty$ as $t \to  + \infty$ and
$$
t= \int \frac {\lambda d \phi}{\phi +O(1)}
$$
which implies that $\phi$ is asymptotic to $e^{\frac t{\lambda}}$
as $t \to + \infty.$ And then from the equation $\phi^{\prime}$ is
asymptotic to $\frac 1{\lambda}e^{\frac t{\lambda}}.$ Since the
metric $g_{\alpha \bar{\beta}}=\partial_{\alpha}
\partial_{\bar{\beta}} u(t)$ is rotationally symmetric, it is clear that
the straight line through the origin are geodesics and the distance
function $d$ from the origin $O$ is a function of $t$ only and is given by
$$
\arraycolsep=1.5pt
\begin{array}{rl}
d= & \displaystyle \int_0^r \sqrt{g_{1 \bar{1}}(r)}dr
\\[4mm]
=& \displaystyle \frac 1 2 \int_{-\infty}^{t}
\sqrt{\phi^{\prime}}dt
\end{array}
$$
by using (7.2). It then follows that the geodesic distance $d$ is asymptotic
to $\sqrt{\lambda}e^{\frac t{2 \lambda}}$ as $t \to +\infty.$ By (7.3),
(7.4) and  the asymptotic behaviors of $\phi$ and $\phi^{\prime}$ as
$t \to +\infty,$ we have the scalar curvature
$$
\arraycolsep=1.5pt
\begin{array}{rl}
\displaystyle R= & \displaystyle g^{\alpha \bar{\beta}}R_{\alpha
\bar{\beta}}
\\[4mm]
= & \displaystyle (e^t \phi^{-1} \delta_{\alpha \beta}+
\bar{z}_{\beta} z_{\alpha}((\phi^{\prime})^{-1}- \phi^{-1}))
\partial_{\alpha}\partial_{\bar{\beta}}(nt-(n-1)  \log
\phi -\log \phi^{\prime})
\\[4mm]
= & \displaystyle O(\phi^{-1})
\\[4mm]
= & \displaystyle O(\frac 1{d^2}).
\end{array}
$$
And by (7.4) and the asymptotic behavior of $\phi$ as $ t \to
+ \infty,$ we have the volume growth
$$
\arraycolsep=1.5pt
\begin{array}{rl}
\displaystyle Vol(B(O,d))= & \displaystyle \omega_{2n-1}
\int_{0}^{d} \det (g_{\alpha \bar{\beta}}) r^{2n-1}dr
\\[4mm]
= & \displaystyle  \frac {\omega_{2n-1}}{2} \int_{-\infty}^t
\phi^{\prime}(\phi)^{n-1}dt
\\[4mm]
= & \displaystyle \frac {\omega_{2n-1}}{2n} \phi^n
\\[4mm]
= & \displaystyle \frac {\omega_{2n-1}}{2n}
(e^{\frac t{\lambda}})^n+ \mbox{ lower order terms}
\\[4mm]
= & \displaystyle \frac {\omega_{2n-1}}{2n \lambda^n}d^{2n}+
\mbox{ lower order terms, }
\qquad \mbox{ as } d \to + \infty,
\end{array}
$$
where $B(O,d)$ is the  geodesic ball  centered at the origin $O$
and of radius $d$ with respect to the K\"{a}hler metric $g_{\alpha
\bar{\beta}}=\partial_{\alpha} \partial_{\bar{\beta}}u(t),$ and
$\omega_{2n-1}$ is the area of the unit sphere $S^{2n-1}.$ The
standard volume comparison theorem thus gives $$Vol(B(x,d)) \geq
\frac {\omega_{2n-1}}{2n\lambda^n}d^{2n}$$ for any $x \in M, d >
0$. \vskip 0.2cm Therefore the complete K\"{a}hler metrics
$g_{\alpha \bar{\beta}}=\partial_{\alpha}
\partial_{\bar{\beta}}u(t)$ obtained by solving (7.8) for every
$\lambda > 1$ satisfy all the assumptions of Theorem 3.

\end{document}